\documentclass{article}
\usepackage[T2A]{fontenc}
\usepackage[cp1251]{inputenc}
\usepackage[english,russian]{babel}
\usepackage[tbtags]{amsmath}
\usepackage{amsfonts,amssymb,mathrsfs,amscd}
\usepackage{graphics}
\usepackage[hyper]{msb-a}
\JournalName{}
\theoremstyle{plain}
\newtheorem{theorem}{Теорема}
\newtheorem{lemma}{Лемма}
\newtheorem{propos}{Предложение}

\newtheorem{corollary}{Следствие}
\theoremstyle{definition}
\newtheorem{definition}{Определение}
\newtheorem{proof}{Доказательство}
\newtheorem{remark}{Замечание}
\newtheorem{comm}{Комментарий}

\def\const{\mathrm{const}}
\def\loc{{$^{\mathrm{loc}}$}}
\def\sloc{{$^{\mathrm{s-loc}}$}}
\def\RR{\mathbb R}
\def\QQ{\mathbb Q}

\def\NN{\mathbb N}

\def\II{\mathbb I}
\def\C{\varkappa}
\def\sgn{\mathrm{sgn}\,}
\begin{document}

\title{Обобщение теоремы Бертрана на поверхности вращения}
\author[O.\,A.~Zagryadsky]{О.\,А.~Загрядский}
\address{МГУ им.~М.\,В.~Ломоносова}
\email{grpcozag@mail.ru}
\author[E.\,A.~Kudryavtseva]{ Е.\,А.~Кудрявцева}
\address{МГУ им.~М.\,В.~Ломоносова}
\email{eakudr@mech.math.msu.su}
\author[D.\,A.~Fedoseev]{Д.\,А.~Федосеев}
\address{МГУ им.~М.\,В.~Ломоносова}
\email{denfedex@yandex.ru}

\date{}
\udk{514.853}

\maketitle

\begin{fulltext}

\begin{abstract}
В работе доказано обобщение теоремы Бертрана на случай абстрактных поверхностей вращения, не имеющих ``экваторов''. Доказан критерий существования на такой поверхности ровно двух центральных потенциалов (с точностью до аддитивной и мультипликативной констант), для которых все ограниченные орбиты замкнуты и имеется ограниченная неособая некруговая орбита.
Доказан критерий существования ровно одного такого потенциала.
Изучены геометрия и классификация соответствующих поверхностей,
с указанием пары (гравитационного и осцилляторного) потенциалов или
единственного (осцилляторного) потенциала. Показано, что на
поверхностях, не относящихся ни к одному из описанных классов,
потенциалов искомого вида не существует.

\end{abstract}


\begin{center} \bf
{D.\,A.~Fedoseev, E.\,A.~Kudryavtseva, O.\,A.~Zagryadsky} \medskip \\
{\large Generalization of Bertrand's theorem to surfaces of revolution} \bigskip
\end{center}

\begin{abstract}
The generalization of Bertrand's theorem to abstract surfaces of 
revolution without ``equators'' is proved. We prove a criterion for the existence on such 
a surface of exactly two central potentials (up to an additive and a multiplicative 
constants) all of whose bounded nonsingular orbits are closed and which admit a bounded 
nonsingular noncircular orbit. A criterion for the existence of a unique such potential is
proved. The geometry and classification of the corresponding surfaces are described, with 
indicating either the pair of (gravi\-ta\-ti\-onal and oscillator) potentials or the unique 
(oscillator) potential. The absence of the required potentials on any surface which does 
not meet the above criteria is proved.



\end{abstract}

\markright{Обобщение теоремы Бертрана на поверхности вращения}

\footnotetext[0]{Работа выполнена в Московском Государственном университете им.\ М.В.\ Ломоносова при поддержке Программы Президента
РФ поддержки ведущих научных школ (грант № НШ-3224.2010.1), Программы
``Развитие научного потенциала высшей школы'' (грант № 2.1.1.3704),
Гранта РФФИ № 10–01–00748-а, ФЦП ``Научные и научно-педагогические кадры инновационной России'' (грант № 14.740.11.0794) и гранта Правительства РФ для господдержки научных исследований, проводимых под руководством ведущих ученых, в ФГБОУ ВПО ``Московский государственный университет имени М.В.~Ломоносова'' по договору № 11.G34.31.0054.}

\section{Введение}

\subsection{Классические результаты}
\
Во второй половине XIX века была поставлена и решена следующая задача о движении точки в евклидовом пространстве $\mathbb{R}^3$ (Ж.Бертран \cite{Bert1}):
{\sl Найти закон силы притяжения, если она зависит только от расстояния и заставляет свою точку приложения описывать замкнутую кривую, каковы бы ни были начальные условия, если только начальная скорость точки меньше некоторого предела}.

Заметим, что, хотя формулировалась эта задача для движения точки в евклидовом пространстве $\mathbb{R}^3,$ в силу центральности рассматриваемых потенциалов ее решение сводится к рассмотрению движения точки по евклидовой плоскости $\mathbb{R}^2.$

Была также поставлена и решена сходная задача (Г.Кёнигс \cite{Koenig1}): {\sl Зная, что сила, вызывающая движение планеты вокруг Солнца, зависит только от расстояния и такова, что она заставляет свою точку приложения описывать алгебраическую кривую, каковы бы ни были начальные условия, причем существуют ограниченные неособые некруговые орбиты, найти закон этой силы}.

Обе задачи были решены в евклидовом пространстве $\mathbb{R}^{3}$ и
ответ одинаков \--- закон притяжения может иметь либо ньютоновский
(т.е.\ гравитационный) вид с силой притяжения $-\frac{G}{r^{2}},$
либо гуковский (т.е.\ осцилляторный) вид с силой притяжения $-kr,$
где $r$ \--- расстояние от точки до притягивающего центра. Другими
словами, уравнения движения точки имеют вид
$\frac{d^2}{dt^2}\vec{r}=-\frac{G}{r^3}\vec{r}$ либо
$\frac{d^2}{dt^2}\vec{r}=-k\vec{r},$ где $\vec{r}=\vec{r}(t) \in
\mathbb{R}^3$ \--- радиус-вектор точки (планеты); $r=|\vec{r}|,
G=\const>0, k=\const>0.$ Первая задача была решена  Ж.Бертраном и
Г.Дарбу \cite{Bert1, Dar2}, см. также \cite{Dar1,Dar3}. Вторая решена Г.Кёнигсом \cite{Koenig1}.

В работе \cite{Bert1} Ж.Бертраном была сформулирована и доказана следующая теорема (в действительности, при дополнительном предположении о том, что центральный потенциал является сильно замыкающим, см.\ определение \ref{closing}).

\begin{theorem} [(Ж.Бертран, 1873 \cite{Bert1})] \label{bert1}
В евклидовом пространстве существуют ровно $2$ закона притяжения с
аналитическим центральным потенциалом, при которых всякая траектория
точки $P,$ движущейся вокруг неподвижной точки $O$ (при условии, что
координаты начального положения точки и компоненты её начальной
скорости не пропорциональны и начальная скорость точки меньше
некоторого предела, зависящего от начального положения точки $P$),
является замкнутой, причем необязательно несамопересекающейся. Этими
законами являются закон Ньютона с силой притяжения
$F_1=-\frac{G}{r^{2}}$ и закон Гука с силой притяжения $F_2=-kr,$ где
$G>0,$ $k>0.$ Для закона сил $F_\beta$ неособые ограниченные
некруговые орбиты задаются периодическими функциями $r=r(\varphi)$ с
минимальным положительным периодом $\Phi_\beta=2\pi/\beta,$
$\beta=1,2.$
\end{theorem}

Причем в случае обоих потенциалов (ньютоновского, т.е.\
гравитационного, $V_1(r)=-\frac{G}{r}$ и гуковского, т.е.\
осцилляторного, $V_2(r)=k\frac{r^2}{2}$) геометрический вид орбит
будет одним и тем же: это будут конические сечения, а в случае
замкнутых орбит --- эллипсы с фокусом или центром в точке притяжения.
Аналог теоремы \ref{bert1} для сферы и плоскости Лобачевского доказал
Г.~Либман~\cite{Lib2}.

Рассуждения, содержащиеся в работах \cite{Bert1,Lib1,Lib2}, основаны
на следующем техническом утверждении, которое вытекает из работы
\cite{Bert1}, и которое мы будем называть технической теоремой Бертрана.

\begin{theorem} [(Ж.Бертран \cite{Bert1}, техническая теорема)] \label{bert2}
Рассмотрим одно\-па\-ра\-мет\-ри\-чес\-кое семейство дифференциальных
уравнений $\frac{d^{2}z}{d\varphi^{2}}+z=\frac{1}{K^2}\Psi(z)$ на
луче $z > 0$ с параметром $K\in \mathbb{R}\setminus\{0\},$ где
$\Psi=\Psi(z)$ \--- аналитическая функция, такая что $\Psi(z)>0.$
(Функция $\frac{1}{K^2}\Psi(z)$ называется силовой или функцией
внешних сил, а функция $-z$ --- центробежной силой или внутренней
силой). Функцию $\Psi(z)$ назовем {\em рационально замыкающей},
если {\rm (i)} для всех $K$ все ограниченные непостоянные решения $z=z(\varphi)$
являются периодическими функциями с периодами, соизмеримыми с $2\pi,$
{\rm (ii)} всякая точка $z>0$ является невырожденным устойчивым положением
равновесия уравнения при $|K|=\sqrt{\Psi(z)/z}.$ Существуют две и
только две рационально замы\-ка\-ю\-щие функции $\Psi$ с точностью
до мультипликативной константы: $\Psi(z)=\frac{A}{z^{\beta^{2}-1}},$
где $\beta\in\{1,2\},$ $A>0$ --- произвольная мультипликативная
константа. При этом все ограниченные непостоянные решения
являются периодическими функциями с минимальным положительным
периодом $\Phi=2\pi/\beta.$
\end{theorem}

В действительности, теорема \ref {bert2} и некоторые ограничения в
ней на функцию $\Psi(z)$ и мультипликативную константу $A$ (а именно:
$\Psi(z)>0,$ (ii), $A>0$) не были явно сформулированы в \cite{Bert1},
однако именно при этих ограничениях эта теорема доказана в \cite{Bert1}.
Теорема \ref{bert2} является в некотором смысле переформулировкой
теоремы \ref{bert1} (и ее аналога для сферы и плоскости Лобачевского,
см.\ \cite {Lib1,Lib2}): здесь $z=\frac{1}{r},$ функция $z(\varphi)$
характеризует зависимость $r(\varphi)$ расстояния от угла,
$2\pi/\beta-$периодичность которой (при рациональном $\beta>0$)
отвечает за замкнутость траектории движения точки, а функция
$-\Psi(z)$ есть производная от потенциала $V(\frac{1}{z})$ по
переменной $z.$

В дальнейшем $1-$параметрическое семейство дифференциальных уравнений
вида $\frac{d^{2}z}{d\varphi^{2}}+z=\frac{1}{K^2}\Psi(z),$ $K \neq
0,$ будем называть {\it семейством уравнений Бертрана}. Из теоремы
\ref{bert2} следует, что если уравнения орбит точки $z(\varphi)$
образуют семейство уравнений Бертрана, где $K$ --- значение интеграла
кинетического момента, то из условия замкнутости ограниченных орбит
(и существования таких орбит) следует, что потенциал $V(r)$ имеет
один из двух определенных видов (с точностью до аддитивной и
мультипликативной констант).

Константу $\beta\in\{1,2\}$ из теорем \ref {bert1}, \ref {bert2}
назовем {\it постоянной Бертрана}, отвечающей данному закону сил. Она
равна 1 для гравитационного и 2 для осцилляторного законов на
следующих классических поверхностях: евклидовой плоскости (см.\
теорему \ref {bert1}), сфере и плоскости Лобачевского (см.\ \cite
{Lib2} или теорему \ref {our1} при $c\ne0,$ $\xi=1$).

\begin{remark} \label{rem:xi}
Если в условии (i) теоремы \ref {bert2} потребовать соизмеримость всех периодов с числом $2\pi\xi$ (вместо $2\pi$) при некотором $\xi>0$, получим определение {\it $2\pi\xi$-замыкающей} функции. Как показывают доказательства технической теоремы \ref {bert2} в работах \cite {Bert1,Lib1}, ее утверждение останется верным, если в условии (i) потребовать лишь попарную соизмеримость периодов (вместо соизмеримости с $2\pi$).
Поэтому при иррациональном $\xi$ не существует ни одной 
$2\pi\xi$-замыкающей функции $\Psi.$ Отсюда следует обобщение
теоремы \ref {bert1} Бертрана на случай всех (не обязательно
рациональных) конусов (см.\ следствие \ref {cone} и \S\ref{sec:cone}), а также обобщение
теоремы Либмана \cite {Lib2} на случай всех (не обязательно
рациональных) ``вещественных разветвленных накрытий'' сферы и
плоскости Лобачевского (см.\ теорему \ref {our1} при $c\ne0$).
\end{remark}

\subsection{Постановка задачи и описание Дарбу всех пар Бертрана}
\
Рассмотрим вместо евклидовой плоскости некоторую {\it абстрактную поверхность вращения}, т.е.\ двумерную поверхность $S \thickapprox (a,b) \times S^1$ с римановой метрикой
\begin{equation}\label{metric}
ds^2=dr^2+f^2(r)d\varphi^2, \quad (r,\varphi \bmod{2\pi})\in(a,b)\times S^1,
\end{equation}
где $f=f(r)$ --- бесконечно гладкая и положительная функция на
интервале $(a,b),$ $-\infty \le a < b \le \infty.$ Рассмотрим систему
на этой поверхности, состоящую из неподвижного притягивающего центра
(``Солнца'') и притягиваемой точки (``планеты''). Возникает
естественный вопрос: что можно сказать о виде центральных потенциалов
$V(r)$ (законов притяжения), для которых все ограниченные орбиты
замкнуты (причем такие орбиты существуют)? Фактически Дарбу \cite
{Dar1,Dar3} и Перлик~\cite {Perl} получили ответ на этот вопрос для
поверхностей без экваторов (см.\ теорему \ref {dar} и замечание \ref
{d=0}), хотя Дарбу поставил его лишь для настоящих (неабстрактных)
поверхностей вращения в $\RR^3.$ Возникают также вопросы об изучении
геометрических и динамических свойств этих поверхностей (частично рассмотренные в~\cite {KozHar,Koz1,BorMam,Sant1,Ball}) и их классификации с точностью до изометричности, подобия и других отношений эквивалентности.

Хочется понять также, на какие абстрактные поверхности вращения
обобщается теорема \ref{bert1} Бертрана: на каких поверхностях
существует ровно два центральных потенциала, которые обеспечивают
замкнутость всех ограниченных траекторий. Существуют ли поверхности,
на которых число требуемых центральных потенциалов (с точностью до
аддитивной и мультипликативной констант) отлично от нуля и двух?
Частичный ответ на этот вопрос недавно был получен М.~Сантопрете 
\cite {Sant1} (см.\ теорему \ref {sant}).

Для определения замыкающего потенциала введем несколько определений.

В дальнейшем под \textit{траекторией} будем понимать решение
$\vec{r}(t)$ уравнения движения, определенное на максимальном по
включению интервале $(t_0,t_1)\subset\mathbb{R}^1,$ под
\textit{орбитой} --- образ этого отображения $O=\{\vec{r}=\vec{r}(t)
\mid t \in (t_0,t_1)\} \subset S,$ под \textit{фазовой траекторией}
--- функцию $(\vec{r}(t), \dot{\vec{r}}(t))$ со значениями в
касательном пространстве $TS,$ а под \textit{фазовой орбитой} --- ее
образ
$\tilde{O}=\{(\vec{r}=\vec{r}(t),\dot{\vec{r}}=\dot{\vec{r}}(t)) \mid
t \in (t_0,t_1)\}\subset TS.$

\begin{definition}\label{traj}
(a) Орбита точки, движущейся по поверхности $S$ по закону сил,
заданному центральным (т.е.\ зависящим только от $r$) потенциалом
$V(r),$ называется \textit{круговой}, если она совпадает с орбитой
действия группы вращений. Траектория называется \textit{круговой},
если соответствующая ей орбита является круговой. Круговую орбиту
назовем \textit{сильно устойчивой}, если на ней функция
$V(r)+\frac{K^2_0}{2f^2(r)},$ называемая \textit{эффективным
потенциалом}, имеет невырожденный локальный минимум при некотором
$K_0\neq0.$

Замкнутую орбиту назовем \textit{орбитально устойчивой}, если
отвечающая ей фазовая траектория орбитально устойчива для ограничения
системы на множество уровня кинетического момента, содержащее эту фазовую 
траекторию.

(b) Траектория называется \textit{ограниченной}, если она определена
на всей оси времени $t \in \mathbb{R}^1,$ и ее образ содержится в
некотором компакте $[r_1,r_2]\times S^1 \subset (a,b) \times S^1.$
Орбита называется \textit{ограниченной}, если соответствующая ей
траектория ограничена.

(c) Траекторию (и соответствующую орбиту) назовем \textit{особой},
если значение интеграла кинетического момента $K$ на этой траектории
равно $0,$ т.е.\ $\varphi = \const.$
\end{definition}

\begin{remark} \label{rem:circle}
(a) Пусть имеется двумерная поверхность $S$ с метрикой (\ref{metric}), где $f(r)$ --- гладкая функция, и центральный гладкий потенциал $V(r)$ на ней. Проверка показывает, что окружность $\{r_0\}\times S^1$ является круговой орбитой тогда и только тогда, когда $r_0$ является критической точкой эффективного потенциала $V(r)+\frac{K_0^2}{2f^2(r)},$ где $K_0$ --- значение интеграла кинетического момента на соответствующей траектории.
Заметим сразу, что $f>0.$ Пусть теперь $f'(r_0)$ и $V'(r_0)$ имеют одинаковый знак или одновременно равны нулю. Получаем, что в этом и только в этом случае окружность $\{r_0\}\times S^1$ является круговой орбитой, причем значение кинетического момента на соответствующей траектории равно $K_0=\pm \sqrt{\frac{V'(r_0)}{f'(r_0)}f^3(r_0)}$ или любое $K_0 \neq 0$ соответственно.

(b) Окружность $\{r\}\times S^1$ является сильно устойчивой круговой
орбитой тогда и только тогда, когда она является круговой орбитой и
соответствующая фазовая орбита является боттовским критическим
множеством локальных минимумов ограничения интеграла энергии на
поверхность уровня интеграла кинетического момента в фазовом
пространстве. Поэтому сильно устойчивые круговые орбиты орбитально устойчивы.
\end{remark}

Заметим, что если траектория неограничена в смысле определения \ref{traj}(b), то
планета ``выходит на край поверхности''.

Мы изучим следующие пять классов центральных потенциалов на
абстрактных поверхностях вращения, а также опишем все поверхности без
экваторов, допускающие такие потенциалы.

\begin{definition} \label{closing}
Пусть $V(r)$ \--- центральный потенциал на поверхности $S$ с
метрикой (\ref{metric}). Будем называть его \textit{замыкающим},
если
\begin{enumerate}
\item[$(\exists)$] существует неособая ограниченная некруговая орбита $\gamma$ в $S$;
\item[$(\forall)$] всякая неособая ограниченная орбита в $S$ замкнута.
\end{enumerate}
Потенциал $V(r)$ будем называть \textit{локально замыкающим}, если
\begin{enumerate}
\item[$(\exists)$\loc] существует сильно устойчивая круговая орбита $\{r_0\}\times S^1$ в $S$;
\item[$(\forall)$\loc] для всякой сильно устойчивой круговой орбиты $\{r_0\}\times S^1$ в $S$
существует $\varepsilon > 0,$ такое что всякая неособая ограниченная
орбита, целиком лежащая в кольце $[r_0-\varepsilon, r_0+\varepsilon]
\times S^1$ и имеющая уровень кинетического момента в интервале
$(K_0-\varepsilon,K_0+\varepsilon),$ является замкнутой, где $K_0$
--- значение кинетического момента на соответствующей круговой
траектории.
\end{enumerate}
Потенциал $V(r)$ будем называть \textit{полулокально замыкающим},
если выполнены условия $(\exists),$ $(\forall)$\loc\ и следующее
условие:
\begin{enumerate}
\item[$(\forall)$\sloc] любая неособая ограниченная орбита в кольце $U = [a',b'] \times S^1$
с уровнем кинетического момента $\hat K$ является замкнутой, где $a'
:= \inf{r|_{\gamma}},$ $b' := \sup{r|_{\gamma}},$ $\gamma$ --
ограниченная орбита из $(\exists),$ $\hat K$ --- значение
кинетического момента на ней.
\end{enumerate}
Потенциал $V(r)$ назовем \textit{сильно} 
 (соответственно \textit{слабо}) \textit{замыкающим}, 
если выполнено условие $(\forall)$\loc\ 
 (соответственно его аналог для всякой орбитально устойчивой круговой орбиты) 
и следующее условие, используемое в работах
\cite {Bert1}, \cite {Dar1,Dar2,Dar3}, \cite {Sant1} 
 (соответственно \cite {Perl}):
любая окружность $\{r\}\times S^1$ является сильно устойчивой 
 (соответственно орбитально устойчивой) 
круговой орбитой.
\end{definition}

\begin{remark} \label{rem*}
(a) Очевидно, что любой замыкающий центральный потенциал $V(r)$
является полулокально замыкающим, а любой сильно замыкающий --- локально 
 и слабо 
замыкающим (см.\ замечание \ref {rem:circle}).
 Нетрудно показывается (при помощи первых интегралов полной энергии и
 кинетического момента), что любой слабо замыкающий центральный
 потенциал является полулокально замыкающим. 
Мы покажем (см.\ теорему \ref {statement}(A) или доказательство теоремы \ref {genbert}, шаг
1), что если $V(r)$ является полулокально замыкающим, то он является
локально замыкающим, причем существует сильно устойчивая круговая
орбита $\{r=r_0\},$ на которой значение интеграла кинетического
момента $K$ совпадает с $K|_{\gamma}$ для неособой ограниченной
некруговой траектории $\gamma$ из $(\exists).$ Тем самым, все
потенциалы из определения \ref {closing} являются локально
замыкающими (без предположения об отсутствии на поверхности экваторов, см.\ ниже).

(b) Мы покажем (теоремы \ref {our1}, \ref {our_1_1}, \ref
{statement}(B)), что если поверхность не имеет {\it экваторов} (т.е.\
таких окружностей $\{r\}\times S^1,$ что $f'(r)=0$),
то все пять понятий замыкающего потенциала $V(r)$ (см.\ определение \ref {closing}) равносильны,
причем тройка $(f(r),V(r),\beta)$ имеет определенный вид (см.\
теорему \ref {dar} или теоремы \ref {our1}, \ref {our_1_1}): либо вид
$(\xi f_{c,0}(\pm(r-r_0)),V_{c,0,i}(\pm(r-r_0)),i\xi)$ для некоторых
$i\in\{1,2\},$ $c\in\RR$ и $\xi \in \QQ\cap\RR_{>0},$ либо вид
$(f_{c,d}(\pm(r-r_0))/\mu,V_{c,d,2}(\pm(r-r_0)),2/\mu)$ для некоторых
$c\in\RR,$ $d\in\RR\setminus\{0\}$ и $\mu\in\QQ\cap\RR_{>0},$ где
$\{\xi f_{c,0}(r)\}$ и $\{f_{c,d}(r)/\mu\}$ --- двупараметрическое и
трехпараметрическое семейства поверхностей (поверхности Бертрана
первого и второго типов, состоящие из одной или двух связных
компонент). Здесь $\{V_{c,0,1}(r)\}$ и $\{V_{c,d,2}(r)\}$ ---
соответствующие семейства замыкающих центральных потенциалов с
точностью до аддитивной и положительной мультипликативной констант,
$\xi,c,d,\mu \in \RR$ --- параметры семейств (см.\ таблицу \ref
{tab:inter}), $\beta=2\pi/\Phi$ --- постоянная Бертрана.

(c) В отличие от теорем \ref{bert1} и \ref {dar} Бертрана и Дарбу (а
также теоремы Перлика~\cite {Perl}), в наших определениях
замыкающего, локально или полулокально замыкающего потенциала не
требуется, чтобы все окружности $\{r\} \times S^1$ являлись сильно
(или орбитально) устойчивыми круговыми орбитами. Лишь в
одном определении требуется существование хотя бы одной такой окружности.
\end{remark}

Из работы Г.Дарбу \cite {Dar1} 1877 г.\ (см.\ также \cite{Dar3,Perl}) вытекает следующее утверждение.

\begin {theorem} [(Г.Дарбу \cite {Dar1,Dar3}, пары поверхность-потенциал Бертрана)] \label {dar}
Пусть на поверхности $S$ с римановой метрикой {\rm(\ref{metric})}
задан центральный потенциал $V=V(r),$ где $f,V$ --- функции класса
$C^\infty,$ не имеющие критических точек.
Потенциал $V$ является сильно замыкающим {\rm(см.\ определение \ref {closing})} в том и только том случае, когда в координатах $(V,\varphi\mod2\pi)$ риманова метрика на $S$ имеет хотя бы один из следующих видов:
 $$
 ds^2=\frac{A\ dV^2}{\beta^2(AV^2-BV+C)^2}+\frac{d\varphi^2}{AV^2-BV+C}
 $$
или
 $$
 ds^2=\frac{A\ dV^2}{\beta^2(-V-K)^3(A/(-V-K)-BV+C)^2}+\frac{d\varphi^2}{A/(-V-K)-BV+C},
 $$
где $A,B,C,K\in\RR,$ $\beta\in\QQ\cap\RR_{>0}$ --- константы, причем
$-2AV+B>0$ или $-A/(V+K)^2+B>0$ соответственно. При этом неособые
ограниченные некруговые орбиты задаются периодическими функциями
$r=r(\varphi)$ с минимальным положительным периодом
$\Phi=2\pi/\beta.$ Поверхности одного вида, отвечающие наборам $(\beta,A,B,C)$ и
$(\beta/\alpha,\alpha^2A,\alpha^2B,\alpha^2C)$ при
$\alpha\in\QQ\cap\RR_{>0},$ локально изометричны друг другу.
\end{theorem}

В действительности, теорема \ref {dar} и указанные в ней формулы для
римановой метрики, как и ограничения на функции $f,V$ и константы
$A,B,K,$ не были явно сформулированы в работах \cite{Dar1,Dar3}
(посвященных в основном реализуемым поверхностям вращения в $\RR^3,$
а не абстрактным поверхностям вращения). Однако именно при этих
ограничениях эта теорема фактически доказана в \cite{Dar1,Dar3}.
(Точнее, указанные выше формулы для римановой метрики следуют из
формулы (15.17) работы \cite{Dar3} путем подстановки в нее решений
(15.10) и (15.11) с учетом соотношения (15.8) и обозначений (15.4), а
остальные утверждения теоремы \ref {dar} следуют из соотношений
(15.3) и (15.19) работы \cite{Dar3}.) В.~Перлик \cite {Perl} обобщил
теорему \ref {dar} Дарбу (в других обозначениях) на более широкий
класс поверхностей и потенциалов: на поверхности класса $C^5$ без
экваторов и слабо замыкающие центральные потенциалы класса $C^5,$
выразив риманову метрику через координаты $(f,\varphi\mod2\pi)$ (см.\
замечание \ref {d=0} о связи параметров Дарбу и Перлика). Хотя в
работе \cite {Perl} результат сформулирован для еще более широкого
класса потенциалов, он строго доказан лишь для сильно замыкающих потенциалов
и не совсем строго --- для слабо замыкающих потенциалов.

\subsection{Развитие вопроса, результат Сантопрете}

Движение точки в центральном поле сил на различных римановых поверхностях неоднократно исследовалось многими учеными (см., например, сборник статей \cite{BorMamSb}). При этом изучалась не только ``прямая задача динамики'' \--- изучение свойств движения по заданному закону, но и ``обратная задача динамики'' \--- поиск законов движения, согласно которым траектории движения будут иметь определенные свойства (например, являться замкнутыми (задача Бертрана), либо алгебраическими кривыми (задача К\"енигса)). Подробная история этих вопросов приведена в \cite{Schep1}. 

Аналог ньютоновской силы как величины, обратной площади сферы радиуса
$r,$ для пространства $H^{3}$ предложил ещё Н.И.Лобачевский
\cite{Lob1} и Я.Больяи \cite{Bol1}. В 1860 г. П.Серре \cite{Serret}
определил аналог гравитационного потенциала на сфере и решил задачу
Кеплера на ней. В 1870 г. Ф.Шеринг написал аналитическое выражение
для потенциала Ньютона на $H^{3}$ \cite{Sche1}. В 1873 г. Р.Липшиц
рассмотрел движение тела в центральном поле на сфере $S^{2}$ со
стандартной метрикой, однако вместо потенциала $-\frac{1}{\tg r}$ он
рассмотрел потенциал $-\frac{1}{\sin r}.$ Он нашел общее решение этой
задачи в эллиптических функциях \cite{Lip1}. В 1885 г. В.Киллинг
обобщил законы Кеплера на сферу $S^{3},$ оснащенную стандартной
метрикой \cite{Kil1}. Подобно Н.И.Лобачевскому и Я.Больяи, он
рассматривал силу притяжения как величину, обратную площади двумерной
сферы радиуса $r$ в $S^{3}.$ В следующем году эти результаты были
также передоказаны К.Нейманом \cite{Neu1}. В работе \cite{Kil1}
В.Киллинг также доказал, что переменные в задаче Кеплера с двумя
притягивающими центрами на сфере $S^{n}$ со стандартной метрикой
разделяются, что влечёт интегрируемость задачи. В 1902 г. Г.Либман
\cite{Lib1} перенёс эти результаты на $H^{3}.$

В 1940-х годах этот вопрос рассматривался в рамках теории
относительности, а именно решалась квантово-механическая
одночастичная спектральная задача для ньютоновского потенциала на
сфере $S^{3}$ Э.~Шрёдингером и А.Ф.~Сти\-вен\-со\-ном, на $H^{3}$ Л.~Инфельдом и
А.~Шильдом. В 1980-х годах центральные потенциалы в рамках теории
относительности на $S^{3},$ $H^{3},$ $S^{n}$ исследовались
Ю.А.~Курочкиным, В.С.~Отчиком, А.А.~Богушем, Г.~Лимоном. В 1994 году
В.В.~Козлов переоткрыл законы Кеплера для пространств постоянной
секционной кривизны \cite{Koz1}. В этом же году он вместе с
Ю.Н.~Фёдоровым установил интегрируемость классического движения одной
частицы по сфере $S^{n}$ в поле, создаваемом гуковскими потенциалами,
расположенными в $2(n+1)$ точках пересечения сферы с координатными
осями.

Что касается обратной задачи динамики, то в 1870 г.\ Ж.~Бертран \cite{Bert1}
решил задачу о нахождении всех сильно замыкающих (см.\ определение \ref {closing}) центральных потенциалов на евклидовой плоскости. Затем Г.Дарбу в 1877 г.\ фактически описал все пары Бертрана
поверхность-потенциал $(f(r),V(r)),$ т.е.\ все абстрактные
поверхности вращения без экваторов с сильно замыкающими центральными
потенциалами на них (см.\ теорему \ref {dar}). В частности, он
построил обобщенные гравитационный и осцилляторный потенциалы на
``рациональных'' конусах и ``рациональных накрытиях'' сферы $S^2,$
построил обобщенный осцилляторный потенциал на некоторой грушевидной
поверхности (отвечающей дуге $\ell_4$ в нашем пространстве
параметров, см.\ рис.~\ref{ris:areas}) и ее ``рациональных
накрытиях'', и доказал замкнутость всех ее геодезических.
В 1902 г. Г.Либман построил обобщение осцилляторного потенциала на
$S^{3}$ и $H^{3}$ \cite{Lib1}. При этом при обобщении ньютоновского
потенциала он требовал, чтобы выполнялся первый закон Кеплера (т.е.\
чтобы орбиты являлись эллипсами с фокусами в центре поля).
Определение осцилляторного потенциала проводилось из требования,
чтобы движение частицы происходило по эллипсам с центром, совпадающим
с центром поля. В 1903 г.\ Г.Либман доказал обобщения теоремы
\ref{bert1} Бертрана на случаи $S^{2}$ и $H^{2}$ \cite{Lib2}. 
Позже обобщение теоремы Бертрана на пространство $S^{n}$ было доказано П.~Хиггсом \cite{Hig1} в 1979 г. В очередной раз обобщение теоремы Бертрана на $S^{3}$ было доказано Я.Е.~Славяновским~\cite{Slaw} в 1980 г., а на $S^{n}$ и $H^{n}$ М.~Икедой и Н.~Катаямой~\cite{IkeKat} в 1982 г. В 1992 году В.В.~Козлов и А.О.~Харин переоткрыли ньютоновский и гуковский потенциалы для пространства $S^{n}$ как полное решение обобщенной задачи Бертрана \cite{KozHar}. В том же году В.~Перлик \cite {Perl} получил описание всех
слабо замыкающих пар Бертрана $(f(r),V(r))$ класса $C^5$ без
экваторов, аналогичное описанию Дарбу (см.\ теорему \ref {dar}),
выразив их в координатах $(f,\varphi\mod2\pi)$ с помощью параметров
(см.\ замечание \ref {d=0}).

В 2008 году М.~Сантопрете \cite{Sant1} обобщил теорему \ref{bert1}
Бертрана на аналитические поверхности вращения с постоянной гауссовой
кривизной и без экваторов, вложенные в $\mathbb{R}^{3}.$ Он также
доказал, что на остальных поверхностях вращения без экваторов
количество сильно замыкающих центральных потенциалов не превосходит
одного, и указал вид этого потенциала,
а также привел необходимое условие на метрику (которое в
действительности является и достаточным) в виде биквадратного
уравнения на постоянную Бертрана для существования ровно одного
потенциала.

\begin{theorem} [(М.~Сантопрете \cite{Sant1})] \label{sant}
Пусть $S\subset \mathbb{R}^3$ --- двумерная поверхность вращения с
координатами $(r,\varphi \bmod{2\pi})\in(a,b)\times S^1$ с
аналитической римановой метрикой {\rm(\ref{metric})}, причем функция
$f$ не имеет критических точек на $(a,b).$ Тогда в классе
аналитических центральных потенциалов на $S$:
\begin{enumerate}
\item существует не более двух (с точностью до аддитивной и положительной мультипликативной констант) сильно замыкающих потенциалов;
\item их ровно два, если и только если $f''f-(f')^{2}=-\xi^2,$ где $\xi$ --- положительная рациональная константа, причем этими потенциалами являются обобщенный гравитационный $V_1(r)$ и обобщенный осцилляторный $V_2(r)$ потенциалы;
\item
если потенциал единствен, то $-f''f+(f')^{2}=:h$ не константа,
и потенциал имеет вид обобщенного осцилляторного потенциала.
\end{enumerate}
Более того, для любого сильно замыкающего центрального потенциала
неособые ограниченные некруговые орбиты задаются периодическими
функциями $r=r(\varphi)$ с одним и тем же минимальным положительным
периодом $\Phi=2\pi/\beta,$ где $\beta$ --- положительная
рациональная константа, зависящая от потенциала и удовлетворяющая
биквадратному тождеству
$\beta^{4}-5(-f''f+f'^{2})\beta^{2}-5ff''f'^{2}+4f''^{2}f^{2}-3f'''f'f^{2}+4f'^{4}=0$;
различным (с точностью до аддитивной и положительной
мультипликативной констант) потенциалам отвечают различные константы;
потенциалу $V_i(r)$ из п.~2 отвечает константа $\beta_i=i\xi,$
$i=1,2.$
\end{theorem}

\begin{remark} \label{rem:4}
Для поверхности вращения с метрикой (\ref{metric}) скалярная кривизна
Римана $R$ вычисляется по формуле $R/2=c:=-f''/f,$ причем $c$ ---
гауссова кривизна поверхности в случае вложимости поверхности в
$\mathbb{R}^3$ (см.\ также следствие \ref {cor:geom}(C)). Поэтому
$h'=-f'''f+f''f'=f^2c',$ где $h:=-f''f+(f')^{2}.$ Иными словами,
выполнение условия 2 теоремы \ref{sant} Сантопрете влечет постоянство
гауссовой кривизны поверхности $S.$ Описание всех абстрактных
поверхностей вращения постоянной скалярной кривизны и без экваторов
получается из леммы \ref {3.3} и замечания \ref {rem:3.3}.
\end{remark}

\begin{remark} \label{rem:sant}
(a) В работе \cite{Sant1} Сантопрете не формулировал условие
отсутствия у функции $f$ критических точек, а вместо сильно
замыкающих потенциалов рассматривал замыкающие потенциалы. Также он
не формулировал п.~3 и последнее утверждение теоремы \ref {sant} в виде отдельных
утверждений. Однако эти утверждения следуют из его работы, и именно
такие условия необходимы для проведения его доказательства. Часть
``если'' пункта~2 теоремы \ref {sant} легко следует из теорем
Бертрана и Либмана \cite {Lib1,Lib2}, а пп.~1, 2 (часть ``только
если'') и 3 легко следуют из последнего утверждения теоремы,
поскольку $h:=-f''f+f'^{2}=(\beta_1^2+\beta_2^2)/5,$ где $\beta_i$
--- корни биквадратного уравнения, данного выше, а если $h$ не равно
константе, то может существовать не более одного постоянного корня
$\beta.$ Отметим, что биквадратное уравнение из теоремы \ref {sant}
имеет вид $\beta^4-5\beta^2h+3ff'h'+4h^2=0$ и превращается в
биквадратное уравнение Тикочинского~\cite{Tik} в случае метрик
постоянной кривизны (см.\ замечание \ref {rem:4}).

(b) Теорема \ref{sant} рассматривает лишь поверхности, вложенные в
$\mathbb{R}^{3}$ как поверх\-ности вращения, что накладывает
определенные условия на функцию $f(r).$ А именно, для вложимости
поверхности с метрикой (\ref{metric}) в $\mathbb{R}^{3}$ как
поверхности вращения необходимо, чтобы $|f'(r)| \le 1$ (и достаточно,
чтобы $|f'(r)|<1$). В силу (\ref {f}) в условиях второго пункта
теоремы \ref{sant} это неравенство равносильно неравенству $|\xi| \le
1,$ если $c\ge0$ (см.\ замечание \ref {rem:4}) и $\inf f = 0$
 (например, когда интервал $(a,b)$ максимален);
неравенству $|\xi| \le \min
\{\frac{1}{\ch(\sqrt{-c}(b-r_0))},\frac{1}{\ch(\sqrt{-c}(a-r_0))}\},$
если $c<0.$
\end{remark}

\subsection{Некоторые нерешенные задачи}
Перечислим некоторые, по-видимому, нерешенные задачи, тесно связанные с изучаемой в настоящей работе.
\begin{enumerate}
\item[1)] Изучить аналог задачи Бертрана на прямом цилиндре и других
абстрактных поверхностях вращения с экваторами (см.\ замечания
\ref{rem*}(b), \ref {rem:7}), на \textit{псевдоримановых поверхностях
вращения} с псевдоримановой метрикой $ds^2=dr^2-f^2(r)d\varphi^2,$ а
также при наличии двух центральных полей сил --- потенциального и
идеального магнитного (см.\ пример в~\cite[\S5]{BorMam}).
\item[2)] Описать максимальные поверхности вращения, вложенные в $\RR^3$ (с экваторами или
без), обладающие замыкающими центральными потенциалами (см.\
замечания \ref {rem:sant}(b), \ref {rem:7} и следствие \ref
{cor:geom}(B)).
\item[3)] (А.Т.~Фоменко) Рассмотрим гладкое риманово многообразие $M$ с отмеченной точкой на нем \--- ``Солнцем''.  Пусть другая точка (назовем ее ``планетой'') движется по $M$ в соответствии с уравнениями Ньютона под воздействием силы притяжения, зависящей
только от расстояния до ``Солнца''. Предположим для простоты, что
многообразие двумерно (хотя не обязательно является поверхностью вращения) и что неособые траектории движения ограничены (см.\ определение \ref{traj}), хотя и не обязательно замкнуты, т.е.\ каждая лежит внутри некоторого
кольца, определяемого начальными данными движения ``планеты''. Планета
движется внутри кольца и в случае общего положения ее траектория
может всюду плотно покрывать это кольцо. Так будет, например, если
соответствующая динамическая система на кокасательном расслоении $T^{*}M$
к $M$ допускает интеграл движения, регулярные поверхности уровня
которого являются двумерными торами. По этим торам движется
интегральная траектория системы. При проекции кокасательного
расслоения $T^{*}M$ на $M$ торы проектируются в кольца.  Задача ставится
так:  найти вид потенциала (``закон силы притяжения''), обеспечивающий
именно такое движение ``планеты'' внутри колец. Для начала было бы
интересно разобраться со случаем трехосного эллипсоида.  Конечно, в
постановке задачи предполагается, что начальные скорости ``планеты''
меньше некоторой величины и ``планета не падает на Солнце''.
Аналогичным образом, естественно, формулируется и многомерная задача.

При решении задачи следует использовать сведения о топологии
слоений на торы Лиувилля и их бифуркациях (``атомах''), возникающих для
интегрируемых гамильтоновых систем, см. \cite{BolKozFom}, \cite{Fom}, \cite{NguFom}, \cite{Fom2}, \cite{KudNikFom}.
\item[4)] Изучить динамику и топологию слоения Лиувилля в классе задач о
движении частицы по абстрактной поверхности вращения
$S\approx(a,b)\times S^1$ с римановой метрикой {\rm(\ref{metric})} в
поле сил, заданном произвольным (не обязательно замыкающим)
центральным потенциалом. В частности, задача интегрируема, имеет
дополнительный первый интеграл --- кинетический момент, поэтому
(например, в случае изоэнергетической невырожденности) замыкание
почти любой орбиты планеты --- кольцо, в котором эффективный
потенциал меньше фиксированного уровня энергии
(см.\ задачу 3). Построить и изучить переменные действия, слоение
Лиувилля в фазовом пространстве (а точнее, его области, в которой
совместные множества уровня интегралов энергии и кинетического
момента компактны и связны, потому являются торами) на торы Лиувилля,
а на базе слоения (называемой бифуркационным комплексом) изучить
возникающую аффинную структуру, построить целочисленные решетки
переменных действия.
\end{enumerate}

Статья имеет следующую структуру. В \S \ref {sec:results}
формулируются основные результаты статьи (теоремы \ref {our1}, \ref
{our_1_1}, \ref {statement}, \ref {genbert}, следствия \ref {cone},
\ref {cor:geom}, \ref {cor:class} и комментарии \ref {com:2}, \ref
{com:1}). В \S \ref {sec:cone} обсуждается связь задач Бертрана на евклидовой плоскости и на конусе.
В \S \ref {sec:technical} доказывается обобщенная техническая теорема
\ref {genbert} Бертрана. В \S \ref {sec:our} доказываются теоремы
\ref {our1}, \ref {our_1_1}, \ref {statement} и следствия \ref
{cor:geom}, \ref {cor:class}.

Авторы благодарны А.Т.~Фоменко за постановку задачи, А.В.~Болсинову, А.В.~Борисову, М.Д.~Малых, С.Ю.~Немировскому, А.А.~Ошемкову и А.В.~Щепетилову за полезные обсуждения.

\section{Формулировка основных результатов} \label {sec:results}

\subsection{Обобщение теоремы Бертрана на абстрактные поверхности вращения без экваторов}

\begin{theorem} [($C^{\infty}-$гладкая теорема для поверхностей Бертрана первого типа)] \label{our1}
Пусть дана гладкая двумерная поверхность $S,$ диффеоморфная
$(a,b)\times S^1,$ снабженная римановой метрикой {\rm(\ref{metric})}
(т.е.\ абстрактная поверхность вращения). Пусть функция $f$
удовлетворяет тождеству $f''f-(f')^{2}=-\xi^{2},$ где $\xi>0$
рационально, т.е.\ $f$ имеет один из следующих видов:
\begin{equation}\label{f}
f(r)=\xi f_{c}(r-r_0) := \left\{ \begin{aligned} 
 \pm \xi (r - r_0), \quad c=0,\\ 
 \frac{\xi}{\sqrt{c}}\sin(\sqrt{c}(r-r_0)), \quad c>0,\\ 
 \pm \frac{\xi}{\sqrt{-c}}\sh(\sqrt{-c}(r-r_0)), \quad c<0,
 \end{aligned} \right.
\end{equation}
где $c$ --- половина скалярной кривизны Римана этой поверхности; в данном случае кривизна постоянна; $2\pi\xi$ --- полный угол в конической точке поверхности (центре поля).
Пусть, далее, функция
$f'(r)$ не имеет нулей на интервале $(a,b).$
Тогда в классе центральных потенциалов на $S$ существуют два и только
два (с точностью до аддитивной и мультипликативной констант)
полулокально замыкающих (соответственно локально замыкающих,
замыкающих, сильно или слабо замыкающих) потенциала $V_1(r),$
$V_2(r).$

Эти замыкающие потенциалы являются обобщенными гравитационным
$V_1(r)$ и осцилляторным $V_2(r),$ т.е.\ имеют вид
$V_i(r)=(-1)^iA|\theta(r)|^{2-i^2}/i+B,$ $i=1,2,$
где $A>0,$ $B$ --- некоторые константы,
$\theta(r)=-\frac{f'(r)}{f(r)}=\pm\sqrt{\frac{\xi^2}{f^2(r)}-c}.$
Соответствующие потенциалу $V_i(r)$ неособые некруговые ограниченные
орбиты задаются периодическими функциями $r=r(\varphi)$ с минимальным
положительным периодом $\Phi_i=\frac{2\pi}{i\xi},$ $i=1,2.$ При этом
на фазовой траектории, отвечающей круговой орбите $\{r\} \times S^1
\subset (a,b) \times S^1,$ значение кинетического момента $K$ равно
 $K_i=\pm\xi\sqrt{\frac{A_i}{|\theta(r)|^{i^2}}},$ $i=1,2.$
Граничная окружность $\{\hat r\}\times S^1,$ на которой достигается
$\inf f(r)$ (т.е.\ $\sup |\Theta(r)|$), является притягивающим
центром поля с замыкающим потенциалом $V_i(r)$ (т.е.\ на ней
достигается $\inf V_i(r)$).
\end{theorem}

Пусть $\xi = \frac{p}{q}.$ Тогда $q-$листное накрытие $\tilde{S}$
всякой такой поверхности $S$ может быть представлено как
разветвленное $p-$листное накрытие (с одной или двумя точками
ветвления) одной из трех ``базисных'' поверхностей: евклидова
плоскость, сфера, плоскость Лобачевского. Здесь накрывающее
пространство $\tilde{S}$ определено следующим свойством: это такое
разветвленное накрытие минимальной степени, которое локально
изометрично накрывает проколотую плоскость (соответственно сферу,
плоскость Лобачевского).

В виде иллюстрации рассмотрим частный случай теоремы \ref{our1}: случай конуса (т.е.\ $c=0$), не обязательно вложимого в трехмерное объемлющее пространство в виде поверхности вращения (т.е.\ не обязательно $\xi\le1$).

\begin{corollary} [($C^\infty$-гладкие замыкающие центральные потенциалы на конусах)]\label{cone}
Пусть задан стандартный конус $S \approx (0,+\infty)\times S^1$ с римановой метрикой
\begin{equation} \label{conemetric}
ds^2=dr^2+\xi^2r^2d\varphi^2, \quad (r,\varphi \bmod{2\pi})\in(0,+\infty)\times S^1,
\end{equation}
где $\xi>0,$ т.е.\ угол при вершине конуса равен $2\pi\xi.$ Если
$\xi$ не является рациональным, то не существует замыкающих
(соответственно локально, полулокально, сильно или слабо замыкающих)
центральных потенциалов на рассматриваемом конусе. Если $\xi$
рационально, то существует два и только два (с точностью до
аддитивной и мультипликативной констант) замы\-ка\-ю\-щих
(соответственно локально, полулокально, сильно или слабо замыкающих)
центральных потенциала на конусе: гравитационный и осцилляторный
по\-тен\-ци\-алы $V_i(r)=(-1)^iAr^{i^2-2}/i+B,$ $i=1,2,$
где $A>0,$ $B$ --- произвольные константы. Соответствующие неособые
некруговые ограниченные орбиты задаются периодическими функциями
$r=r(\varphi)$ с минимальным поло\-жи\-тель\-ным периодом
$\Phi_i=2\pi/(i\xi).$
Вершина конуса $\{0\}\times S^1$ является притя\-ги\-ва\-ю\-щим
центром поля для каждого потенциала $V_1,V_2.$
\end{corollary}

В частном случае сильно замыкающих аналитических потенциалов теорема
\ref {our1} и следствие \ref {cone} вытекают из усиления технической
теоремы \ref {bert2} Бертрана (см.\ замечание \ref {rem:xi}), а при
дополнительном условии $\xi\le1$ --- из теоремы \ref {sant}
Сантопрете.

\begin{theorem} [($C^{\infty}-$гладкая теорема для поверхностей Бертрана второго типа)] \label{our_1_1}
Пусть дана гладкая двумерная поверхность $S,$ диффеоморфная
$(a,b)\times S^1,$ снабженная римановой метрикой {\rm(\ref{metric})}.
Пусть функция $f$ не удовлетворяет тождеству $f''f-(f')^{2}=-\xi^{2}$
ни для какого рационального $\xi > 0,$ и пусть функция $f(r)$ не
имеет критических точек на $(a,b).$ Тогда существует не более одного
полулокально замыкающего (соответственно локально замыкающего,
замыкающего, сильно или слабо замыкающего) центрального потенциала (с
точностью до аддитивной и мультипликативной констант). При этом
потенциал ровно один (с точностью до аддитивной и мультипликативной
констант) тогда и только тогда, когда 
существует гладкая функция $\theta = \theta (r)$ без нулей на
$(a,b),$ такая что $\theta'(r)>0$ и риманова метрика в координатах
$(\theta,\varphi\mod2\pi)$ имеет вид
 \begin{equation} \label {eq:sol}
 ds^2
 =\frac{d\theta^2}{(\theta^2+c-d\theta^{-2})^2}
 +\frac{d\varphi^2}{\mu^2(\theta^2+c-d\theta^{-2})},
 \end{equation}
где $\mu$ --- положительная рациональная константа, $d$ --- ненулевая
константа, а $c$ --- произвольная вещественная константа.

При этом функция $\theta(r)$ и тройка чисел $(\mu,c,d)$ единственны
(если существуют), и замыкающий потенциал будет являться обобщенным
осцилляторным, т.е.\ иметь вид
 $V_2(r) = \frac{A}{2\theta^2(r)}+B,$ где $A,B\in\RR,$ $A(\theta^4(r)+d)>0.$
Соответствующие неособые некруговые ограниченные орбиты задаются
периодическими функциями $r=r(\varphi)$ с минимальным положительным
периодом $\Phi=\pi\mu.$ На фазовой траектории, отвечающей круговой
орбите $\{r\} \times S^1 \subset (a,b) \times S^1,$ значение
кинетического момента $K$ равно
 $K=\pm\frac1\mu\sqrt{\frac{A}{\theta^4(r)+d}}$;
граничная окружность $\{r_0\}\times S^1,$ на которой достигается
$\inf f(r)$ (т.е.\ $\sup A|\theta(r)|$), является притягивающим
центром поля с замыкающим потенциалом $V_2$ (т.е.\ на ней достигается
$\inf V_2(r)$).
\end{theorem}

\begin{remark} \label {rem:7}
Отметим, что каждая из теорем \ref{our1}, \ref{our_1_1} и \ref
{statement} предполагает лишь $C^\infty$-гладкость (а не обязательно
аналитичность) функций $f(r),V(r),$ и включает в себя пять
утверждений:
\begin{enumerate}
\item [1)] описание замыкающих центральных потенциалов,
\item [2)] описание полулокально замыкающих центральных потенциалов,
\item [3)] описание локально замыкающих центральных потенциалов,
\item [4)] описание сильно замыкающих центральных потенциалов,
\item [5)] описание слабо замыкающих центральных потенциалов.
\end{enumerate}
В условии теорем \ref{our1}, \ref{our_1_1} и \ref {statement} все
пять классов потенциалов совпадают для абстрактных поверхностей
вращения без экваторов (см.\ замечание \ref {rem*}(b)). В теоремах
\ref{our1}, \ref{our_1_1} и \ref {statement} данной работы (как и в
\cite{Bert1}, \cite {Dar1,Dar2,Dar3}, \cite {Perl}, \cite {Sant1} и
других известных нам работах) не рассматриваются поверхности с
экваторами, например, цилиндр. Предположение об отсутствии экваторов
существенно, так как все поверхности вращения с замкнутыми
геодезическими (т.е.\ поверхности Таннери~\cite[теорема 4.13]{Besse})
обладают замыкающим центральным потенциалом, равным константе.
Можно показать, что при отсутствии экваторов условие (\ref {eq:sol})
при $c,d\in\RR,$ как условие на функцию $f,$ равносильно
биквадратному уравнению из теоремы \ref {sant} Сантопрете на
постоянную Бертрана $\beta:=2\pi/\Phi=i/\mu,$ зависящую от типа
сильно замыкающего потенциала $V_i,$ $i=1,2.$
\end{remark}

\begin{remark}\label{d=0}
Если в теореме \ref{our_1_1} положить
$d=0,$
$\mu = \frac{1}{\xi},$ то соответствующие поверхности в
действительности совпадут с поверхностями, описанными в теореме
\ref{our1}, а потенциал --- с осцилляторным потенциалом $V_2(r).$
Пары поверхность-потенциал из теорем \ref {our1} и \ref {our_1_1},
очевидно, совпадают с парами из теоремы \ref {dar} Дарбу.
Трехпараметрическое семейство пар $(f_{\beta,K}(r),V_{K,G}(r))$ из
работы Перлика \cite {Perl} совпадает с описанным в теореме \ref
{our1} семейством для гравитационного потенциала с мультипликативной
константой $A_1=\xi/2,$ где
$\beta:=\xi,$ $K:=-c\mu^2,$ $G:=-2B$ -- параметры семейства в \cite
{Perl}. Четырехпараметрическое семейство пар
$(f_{\beta,D,K}(r),V_{D,K,G}(r))$ из работы Перлика \cite {Perl}
совпадает с описанным в теореме \ref {our_1_1} семейством (при любом
 $d\in\RR$)
для мультипликативной константы
 $A_2=\pm1/(2\mu^2),$
где $\beta:=2/\mu,$ $D:=\mu^2c,$ $K:=-4\mu^4d,$ $G:=-2B$ ---
параметры семейства в \cite {Perl}.
\end{remark}

\subsection {Геометрия и классификация поверхностей Бертрана}

Ниже мы построим единое семейство максимальных аналитических (вообще
говоря, несвязных) поверхностей, содержащих в себе любую из
поверхностей, описанных в теоремах \ref{our1} и \ref{our_1_1}, и изучим геометрию и классификацию этих поверхностей. 

Рассмотрим семейство функций
$Q_{c,d}(\theta):=\theta^2+c-d\theta^{-2},$ $c,d \in \RR.$ Обозначим
 $$
 I_{c,d}:=\{\theta\in\RR\mid\theta<0,\ Q_{c,d}(\theta)>0,\
 Q_{c,d}'(\theta)\ne0 \} =:\bigcup_{k=1}^{k_{c,d}} I_{c,d,k},
 $$
где $I_{c,d,k}=(\theta_{c,d,k, \, min}, \theta_{c,d,k, \, max})
\subset I_{c,d}$ --- максимальный по включению интервал, $k$ -- номер
интервала, $k_{c,d}$ -- количество интервалов, т.е.\
$k_{c,d}:=1$ при $d \ge 0$ и $k_{c,d}:=2$ при $d<0.$ Пусть
$r_{c,d}(\theta)=\int\frac{d\theta}{Q_{c,d}(\theta)}$ --- некоторая
первообразная в области $I_{c,d}.$

\begin{figure}[tbp]
\setlength{\unitlength}{10pt}
\begin{center}
%
\begin{picture}(32,16)(-16,-9)
\thicklines \sl
\put(0,0){\put(-.386,-.223){$\bullet$}} \put(-1.2,-1){$(0,0)$}
\qbezier[800](-16,0)(0,0)(16,0)
     \put(-16,.3){$\ell_2$ \quad Плоскости Лобачевского}
     \put(4,.3){$d=0$} \put(8.7,.3){Полусферы \quad $\ell_1$}
\qbezier[400](-10,-8)(-5,0)(0,0) \put(-7,-5){$\ell_3$} \put(-9.5,-8){$c^2+4d=0$}
\qbezier[40](10,-8)(5,0)(0,0)
\put(8.5,-7.5){$\ell_4$}
\put(-3.5,5.3){Семейство полубесконечных}
\put(-2.6,4.1){поверхностей Бертрана}
\put(-8,4){$\Omega_2$}
\put(-.2,-7.8){$\Omega_4$}
\put(.4,-3.9){Семейство грушевидных}
\put(.7,-5.1){поверхностей Бертрана}
\put(13.5,-6){$\Omega_1$}
\put(-15,-1.7){Два семейства}
\put(-15.6,-2.9){полубесконечных}
\put(-14.8,-4.1){поверхностей}
\put(-14,-5.3){Бертрана}
\put(-15,-7){$\Omega_3$}
\end{picture}
\end{center}
\setlength{\unitlength}{1pt}
\caption{Плоскость параметров $\RR^2(c,d),$ разбитая на подмножества:
области, кривые и точку (на дуге $\ell_4$ нет бифуркации семейства).}
\label{ris:areas}
\end{figure}

\begin{remark}
Функции $r_{c,d}(\theta)$ и максимальные интервалы $I_{c,d,k}$ в
зависимости от $(c,d)\in\RR^2$ могут быть заданы явными формулами. А
именно, обозначим $\Delta=\Delta(c,d):=c^2+4d.$ Касающиеся друг друга
прямая $\{d=0\}$ и парабола $\{\Delta=0\}$ разбивают плоскость
$\RR^2$ на следующие подмножества: области
$\Omega_1=\{\Delta>0,d<0<c\},$ $\Omega_2=\{d>0\},$
$\Omega_3=\{c<0,d<0<\Delta\},$ $\Omega_4=\{\Delta<0\},$ кривые
$\ell_1=\{d=0<c\},$ $\ell_2=\{c<0=d\},$ $\ell_3=\{c<0=\Delta\},$
$\ell_4=\{\Delta=0<c\}$ и точку $\{(0,0)\}$ (см.\
рис.~\ref{ris:areas}). В области $\overline{\Omega_i}$ через
$x_i=x_i(c,d)$ и $y_i=y_i(c,d),$ $(c,d)\in \overline{\Omega_i}, \,
i\in\{1,2,3,4\},$ обозначим следующие функции:
$$
x_1:=\sqrt{\frac{c-\sqrt{\Delta}}{2}}, \qquad x_2=x_3:=\sqrt{\frac{-c+\sqrt{\Delta}}{2}}, \qquad x_4:=\sqrt{\frac{\sqrt{-d}}{2}-\frac{c}{4}},
$$
$$
y_1=y_2:=\sqrt{\frac{c+\sqrt{\Delta}}{2}}, \qquad y_3:=\sqrt{\frac{-c-\sqrt{\Delta}}{2}}, \qquad y_4:=\sqrt{\frac{\sqrt{-d}}{2}+\frac{c}{4}}.
$$
Тогда функция $r_{c,d}(\theta)$ может быть задана следующими формулами:
$$
r_{c,d}(\theta)=\left\{
\begin{aligned}
-\frac{1}{\theta}, \quad (c,d)=(0,0), \\
\frac{1}{y_1^2-x_1^2}\left(-x_1\arctg\frac{\theta}{x_1}+y_1\arctg\frac{\theta}{y_1}\right), \quad (c,d) \in \Omega_1, \\
\frac{1}{y_1}\arctg\frac{\theta}{y_1}, \quad (c,d) \in \ell_1, \\
\frac{1}{x^2_2+y^2_2}\left(\frac{x_2}{2}\ln\left|\frac{\theta-x_2}{\theta+x_2}\right|+y_2\arctg\frac{\theta}{y_2}\right), \quad (c,d)\in\Omega_2, \\
\frac{1}{2x_2}\ln\left|\frac{\theta-x_2}{\theta+x_2}\right|, \quad (c,d)\in \ell_2, \\
\frac{1}{x^2_3-y^2_3}\left(\frac{x_3}{2}\ln\left|\frac{\theta-x_3}{\theta+x_3}\right|-\frac{y_3}{2}\ln\left|\frac{\theta-y_3}{\theta+y_3}\right|\right), \quad (c,d) \in \Omega_3, \\
-\frac{1}{2}\frac{\theta}{\theta^2-x_3^2}-\frac{1}{4x_3}\ln\left|\frac{\theta+x_3}{\theta-x_3}\right|, \quad (c,d) \in \ell_3, \\
\frac{1}{4y_4}\left(\arctg\frac{\theta+x_4}{y_4}+\arctg\frac{\theta-x_4}{y_4}-\frac{1}{2}\ln\frac{(\theta+x_4)^2+y_4^2}{(\theta-x_4)^2+y_4^2}\right), \quad (c,d) \in \Omega_4, \\
\frac{1}{2y_4}\arctg\frac{\theta}{y_4}-\frac{1}{2}\frac{\theta}{\theta^2+y_4^2}, \quad (c,d) \in \ell_4.
\end{aligned}
\right.
$$

Интервалы $I_{c,d,k},$ на которых задана функция $r_{c,d}(\theta),$
имеют вид
$$
I_{c,d}=\bigcup_{k=1}^{k_{c,d}}I_{c,d,k}:=\left\{ \begin{array}{ll}
(-\infty,0), & (c,d)\in \{(0,0)\}\cup \ell_1;\\
(-\infty,-x_3), & (c,d)\in \Omega_2\cup \ell_2;\\
(-\infty,-x_3)\cup(-y_3,0), & (c,d)\in \Omega_3\cup \ell_3;\\
(-\infty,-\sqrt[4]{-d})\cup(-\sqrt[4]{-d},0), & (c,d)\in \ell_3\cup\Omega_4\cup \ell_4\cup\Omega_1.
\end{array}
\right.
$$
Зависимость этих интервалов от $(c,d)\in\RR^2$ и $k\in\{1,k_{c,d}\}$
показана на рис.~\ref{ris:max} и в таблице \ref {tab:inter}, вместе с
множеством значений функции $r_{c,d}|_{I_{c,d,k}}.$
\end{remark}

По построению каждая из функций $r_{c,d}|_{I_{c,d,k}}:I_{c,d,k}\to
r_{c,d}(I_{c,d,k})$ строго монотонна, а потому имеет обратную
функцию, которую обозначим через $\theta_{c,d,k}=\theta_{c,d,k}(r).$
Итак, функция $\theta_{c,d,k}:r_{c,d}(I_{c,d,k})\to I_{c,d,k}$ такова, что
$\theta_{c,d,k}(r_{c,d}(\theta))\equiv\theta.$
Рассмотрим трехпараметрическое семейство абстрактных поверхностей
вращения $S_{c,d}=\cup_{k=1}^{k_{c,d}}S_{c,d,k}$ с римановыми
метриками $ds^2_{\mu,c,d},$ где
\begin{equation}\label{scdk}
\begin{aligned}
S_{c,d,k} &= r_{c,d}(I_{c,d,k}) \times S^1, \quad
ds^2_{\mu,c,d}|_{S_{c,d,k}}=dr^2+\frac{1}{\mu^2}f^2_{c,d,k}(r)d\varphi^2,\\
S_{c,d} &\approx I_{c,d} \times S^1, \quad ds^2_{\mu,c,d}
 =       \frac{d\theta^2}{(\theta^2+c-d\theta^{-2})^2}
 + \frac{d\varphi^2}{\mu^2(\theta^2+c-d\theta^{-2})},
\end{aligned}
\end{equation}
где $c,d\in\RR,$ $\mu>0$ --- параметры, $k\in\{1,k_{c,d}\},$
$f_{c,d,k}(r):=\frac{1}{\sqrt{Q_{c,d}(\theta_{c,d,k}(r))}}.$ Эта
поверхность состоит из $k_{c,d}$ связных компонент. Первую связную
компоненту $(k=1)$ назовем {\it основной}, вторую ($k=2$ при $d<0$)
--- {\it дополнительной}, а семейство (\ref {scdk}) --- {\em
семейством (максимальных) поверхностей Бертрана}.

\begin{figure}[tbp]
\setlength{\unitlength}{10pt}
\begin{center}
%
\begin{picture}(32,10)(-9.5,0)
\put(-10,8){\line(1,0){10}}  
\put(-10,0){\vector(1,0){13}}\put(3.5,-.5){$\theta$}  
\put(0,0){\vector(0,1){10}} \put(.3,9.5){$\phi\mod2\pi$} 
\put(-.8,-.8){$0$} 
\put(.3,-.5){$\ell_1$}
\put(.3,1.8){$\ell_2$}
\put(.3,3.8){$\ell_3$}
\put(.3,5.8){$\ell_4$}
\put(.3,8){$2\pi$}
\put(1,.7){$\Omega_2$}
\put(1,2.7){$\Omega_3$}
\put(1,4.7){$\Omega_4$}
\put(1,6.7){$\Omega_1$}
\thicklines
\qbezier[800](0,0)(-12,4)(0,8) \put(-5,1.2){$\gamma_2$} \put(-4.2,7){$\gamma_1$}
\qbezier[400](0,2)(-5,2.7)(-6,4) \put(-2.5,2.8){$\gamma_3$}
\put(-2.8,4.9){$\II_2$}
\put(-9,4){$\II_1$}
 \put(1,0){
\put(5,7){$\gamma_1=\{\theta=-\sqrt[4]{-d(\phi)}\}$ \qquad \ ({\sl экваторы})}
\put(5,4){$\gamma_2=\{\theta=-x_3(c(\phi),d(\phi))\}$}
\put(5,2.4){$\gamma_3=\{\theta=-y_3(c(\phi),d(\phi))\}$} \put(17,3.2){ \ ({\sl абсолюты})}
 }
\thinlines   
\qbezier[1](-5.8,3.6)(-5.8,3.6)(-5.8,3.6)
\qbezier[1](-5.6,3.4)(-5.7,3.4)(-5.8,3.4)
\qbezier[3](-5,3.2)(-5.3,3.2)(-5.6,3.2)
\qbezier[5](-4.6,3)(-5.1,3)(-5.6,3)
\qbezier[6](-4,2.8)(-4.6,2.8)(-5.2,2.8)
\qbezier[9](-3.2,2.6)(-4.1,2.6)(-5,2.6)
\qbezier[10](-2.8,2.4)(-3.8,2.4)(-4.8,2.4)
\qbezier[13](-1.8,2.2)(-3.1,2.2)(-4.4,2.2)
\qbezier[20](0,2)(-2,2)(-4,2)
\qbezier[18](0,1.8)(-1.8,1.8)(-3.6,1.8)
\qbezier[19](0,1.6)(-1.9,1.6)(-3.8,1.6)
\qbezier[16](0,1.4)(-1.6,1.4)(-3.2,1.4)
\qbezier[13](0,1.2)(-1.3,1.2)(-2.6,1.2)
 \qbezier[10](0,1.0)(-1,1.0)(-2,1.0)
 \qbezier[6](0,.8)(-.6,.8)(-1.2,.8)
 \qbezier[4](0,.6)(-.4,.6)(-.8,.6)
 \qbezier[3](0,.4)(-.3,.4)(-.6,.4)
 \qbezier[2](0,.2)(-.2,.2)(-.4,.2)
\put(0,0) {\put(-.15,-.11){\tiny$\bullet$}}
\put(0,2) {\put(-.15,-.11){\tiny$\bullet$}}
\put(0,4) {\put(-.15,-.11){\tiny$\bullet$}}
\put(0,6) {\put(-.15,-.11){\tiny$\bullet$}}
\put(0,8) {\put(-.15,-.11){\tiny$\bullet$}}
\end{picture}
\end{center}
\setlength{\unitlength}{1pt}
\caption{Объединения интервалов $\II_k:=\cup_\phi I_{c(\phi),d(\phi),k}\times\{\phi\}$
 для семейств поверхностей Бертрана $S_{c(\phi),d(\phi),k}$ ($k=1,2$), где
$(c(\phi),d(\phi)):=8(\sqrt2\cos(\phi+\pi/4),\sqrt2\sin(\phi+\pi/4)-1).$}
\label{ris:max}
\end{figure}
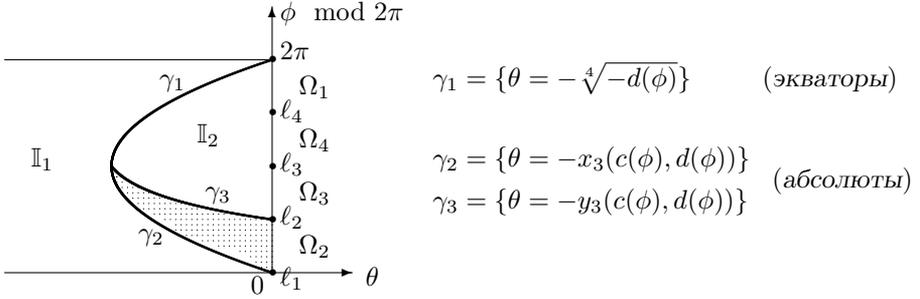

Пусть $\gamma=\gamma(\phi)=(c_\phi,d_\phi),$ $0\le\phi\le2\pi$ ---
некоторая, не обязательно непрерывная, кривая в $\mathbb{R}^2,$
пересекающая всякое множество вида $\{(\lambda c,\lambda^2d)\mid
\lambda>0\}\subset\RR^2$ в единственной точке. Например, образ кривой
$\gamma=\gamma(\phi)$ может являться объединением двух прямых
$\RR\times\{\pm1\}$ и трех точек $\{(-1,0), (0,0), (1,0)\}.$

\begin{theorem} \label{statement} Пусть дана двумерная поверхность $S,$ диффеоморфная
$(a,b)\times S^1,$ снабженная $C^\infty$-гладкой римановой метрикой
{\rm(\ref{metric})}.

{\rm(A)} Любой $C^\infty$-гладкий замыкающий центральный потенциал является полулокально замыкающим, любой сильно замыкающий --- слабо замыкающим, любой слабо замыкающий ---
полулокально замыкающим, а любой полулокально замыкающий --- локально замыкающим.

{\rm(B)} Пусть $f(r)$ не имеет критических точек на $(a,b).$ Тогда
любой $C^\infty$-гладкий локально замыкающий центральный потенциал
является замыкающим и сильно замыкающим, и следующие утверждения
равносильны:
\begin{enumerate}
\item [\rm(a)] на поверхности $S$ имеется замыкающий центральный потенциал $V(r);$
\item [\rm(b)] существует тройка $(\mu,c,d)\in\RR^3$ и функция $\theta=\theta(r)$
без нулей на интервале $(a,b),$ такие что
$\mu\in\mathbb{Q}\cap\mathbb{R}_{>0}$ и выполнено {\rm(\ref{eq:sol})}, т.е.\
\begin{equation} \label{fthetai}
\begin{aligned}
f(r)=\frac{1}{\mu}f_{c,d,k}(\eta(r-r_0)), \ \ \theta(r)&=\eta\theta_{c,d,k}(\eta(r-r_0)), \quad r\in(a,b), \\
(\eta(a-r_0),\eta(b-r_0)) &\subseteq r_{c,d}(I_{c,d,k})
\end{aligned}
\end{equation}
для некоторых $r_0\in\RR,$ $\eta\in\{+,-\}$ и $k\in \{1,k_{c,d}\},$ т.е.\
поверхность $(S,ds^2)$ изометрично и $S^1$-эквивариантно вкладывается
в поверхность Бертрана {\rm(\ref {scdk})} посредством отображения $(r,\varphi)\mapsto(\eta(r-r_0),\varphi)$;
\item [\rm(c)] существует набор $(\mu,\lambda,\phi,r_0,k),$ такой что
$\mu\in\mathbb{Q}\cap\mathbb{R}_{>0},$ $\lambda\in\RR\setminus\{0\},$
$0\le\phi\le2\pi,$ $r_0\in\mathbb{R},$ $k\in\{1,k_{\gamma(\phi)}\},$
$|\lambda|=1$ при $\gamma(\phi)=(0,0),$ причем функция $f(r)$ имеет
вид:
$$
 f(r)=\frac{1}{\mu|\lambda|}f_{\gamma(\phi),k}(\lambda(r-r_0)), \quad
 (\lambda(a-r_0),\lambda(b-r_0)) \subseteq r_{\gamma(\phi)}(I_{\gamma(\phi),k}),
$$
т.е.\ поверхность $(S,\lambda^2ds^2),$ подобная исходной с
коэффициентом подобия $\lambda,$ изометрично и $S^1$-эквивариантно
вкладывается в поверхность
$(S_{\gamma(\phi),k},ds^2_{\mu,\gamma(\phi),k})$ посредством
отображения $(r,\varphi)\mapsto(\lambda(r-r_0),\varphi).$
\end{enumerate}
Наборы чисел в {\rm(b)} и {\rm(c)} единственны. Наборы
$(\mu,\eta,r_0,k)$ из {\rm(b)} и $(\mu,\lambda/|\lambda|,r_0,k)$ из
{\rm(c)} совпадают. Параметры $c,d$ и функция $\theta(r)$ из {\rm(b)}
связаны с параметрами $\lambda,\phi$ из {\rm(c)} соотношениями
$c=\lambda^2c_\phi,$ $d=\lambda^4d_\phi,$
$\theta(r)=\lambda\theta_{\gamma(\phi),k}(\lambda(r-r_0)).$

Потенциал $V=V(r)$ из {\rm(a)} имеет вид
 $V_i(r)=(-1)^iA_i|\theta(r)|^{2-i^2}/i+B_i$
где $i\in\{1,2\}$ при $d=0,$ $i=2$ при $d\ne0,$ $A_i,B_i\in\RR$ ---
любые мультипликативная и аддитивная константы, такие что $A_1>0,$
$A_2(\theta^4(r)+d)>0.$

{\rm(C)} Отвечающие замыкающему центральному потенциалу $V_i$
неособые некруговые ограниченные орбиты задаются периодическими
функциями $r=r(\varphi)$ с минимальным положительным периодом
$\Phi_i=2\pi\mu/i,$ где $\mu$ то же, что и в {\rm(b)} и {\rm(c)},
$i=1,2.$ При этом на фазовой траектории, отвечающей круговой орбите
$\{r\} \times S^1 \subset (a,b) \times S^1,$ значение кинетического
момента $K$ равно
 $K_i=\pm\frac1\mu\sqrt{\frac{A_i}{|\theta(r)|^{i^2}+d}}$;
граничная окружность $\{\hat{r}\}\times S^1,$ на которой достигается
$\inf f(r)$ (т.е.\ $\sup A_i|\theta(r)|$), является притягивающим
центром поля с потенциалом $V_i$ (т.е.\ на ней достигается $\inf
V_i(r)$).
\end{theorem}

\begin{table}
\begin{center} 
\rotatebox{90}{\begin{tabular}{|c|c|c|c|c|c|c|c|c|} \hline
\multicolumn{3}{|c|}{$(c,d)\in \ell_2$} &
\multicolumn{3}{|c|}{$(c,d)\in \Omega_2$} &
\multicolumn{3}{|c|}{$(c,d)\in \ell_1$}
\\
\multicolumn{3}{|c|}{Плоскость Лобачевского} &
\multicolumn{3}{|c|}{Семейство полубесконечных} &
\multicolumn{3}{|c|}{Полусфера
}
\\
\multicolumn{3}{|c|}{
} & \multicolumn{3}{|c|}{поверхностей Бертрана} &
\multicolumn{3}{|c|}{}
\\ \hline
$I_{c,d}$ & \multicolumn{2}{|c|}{$(-\infty,-x_2)$} & $I_{c,d}$ &
\multicolumn{2}{|c|}{$(-\infty,-x_2)$} & $I_{c,d}$ &
\multicolumn{2}{|c|}{$(-\infty,0)$}
\\
$r_{c,d}(I_{c,d})$ & \multicolumn{2}{|c|}{$(0, \infty)$} &
$r_{c,d}(I_{c,d})$ & \multicolumn{2}{|c|}{$(-\frac{\pi
y_2}{2\sqrt{\Delta}}, \infty)$} & $r_{c,d}(I_{c,d})$ &
\multicolumn{2}{|c|}{$(-\frac{\pi}{2y_1},0)$}
\\
$F_{c,d}(I_{c,d})$ & \multicolumn{2}{|c|}{$(0, \infty)$} &
$F_{c,d}(I_{c,d})$ & \multicolumn{2}{|c|}{$(0, \infty)$} &
$F_{c,d}(I_{c,d})$ & \multicolumn{2}{|c|}{$(0, \frac{1}{y_1})$}
\\
$R_{c,d}(I_{c,d})$ & \multicolumn{2}{|c|}{$2c$} & $R_{c,d}(I_{c,d})$
& \multicolumn{2}{|c|}{$(2c, -\frac{2\Delta}{x_2^2})$} &
$R_{c,d}(I_{c,d})$ & \multicolumn{2}{|c|}{$2c$}
\\
$\sgn A_2$ & \multicolumn{2}{|c|}{$+$} & $\sgn A$ &
\multicolumn{2}{|c|}{$+$} & $\sgn A_2$ & \multicolumn{2}{|c|}{$+$}
\\ \hline
\multicolumn{3}{|c|}{$y_2=y_3=0,$ $x_2=x_3=\sqrt{-c}$} &
\multicolumn{3}{|c|}{$x_2>0,$ $y_2>0$} &
\multicolumn{3}{|c|}{$x_1=x_2=0,$ $y_1=y_2=\sqrt{c}$}
\\ \cline{1-9}
\multicolumn{3}{|c|}{$(c,d)\in \Omega_3$} &
\multicolumn{3}{|c|}{$(c,d)=(0,0)$} &
\multicolumn{3}{|c|}{$(c,d)\in \Omega_1$} \\
\multicolumn{3}{|c|}{Два семейства полубесконечных} &
\multicolumn{3}{|c|}{Евклидова плоскость} &
\multicolumn{3}{|c|}{Семейство грушевидных}
\\
\multicolumn{3}{|c|}{поверхностей Бертрана} & \multicolumn{3}{|c|}{}
& \multicolumn{3}{|c|}{поверхностей Бертрана}
\\ \hline
$I_{c,d,k}$ & $(-\infty,-x_3)$ & $(-y_3,0)$ & $I_{c,d}$ &
\multicolumn{2}{|c|}{$(-\infty,0)$} & $I_{c,d,k}$ &
$(-\infty,-\sqrt[4]{-d})$ & $(-\sqrt[4]{-d},0)$
\\
$r_{c,d}(I_{c,d,k})$ & $(0, \infty)$ & $(-\infty,0)$ &
$r_{c,d}(I_{c,d})$ & \multicolumn{2}{|c|}{$(0, \infty)$} &
$r_{c,d}(I_{c,d,k})$ & $(\frac{-\pi}{2(x_1+y_1)},-r_1)$ & $(-r_1,0)$
\\
$F_{c,d}(I_{c,d,k})$ & $(0, \infty)$ & $(\infty,0)$ &
$F_{c,d}(I_{c,d})$ & \multicolumn{2}{|c|}{$(0, \infty)$} &
$F_{c,d}(I_{c,d,k})$ & $(0, \frac{1}{x_1+y_1})$ &
$(\frac{1}{x_1+y_1},0)$
\\
$R_{c,d}(I_{c,d,k})$ & $(2c,-\frac{2\Delta}{x_3^2})$ &
$(-\frac{2\Delta}{y_3^2},\infty)$ & $R_{c,d}(I_{c,d})$ &
\multicolumn{2}{|c|}{$0$} & $R_{c,d}(I_{c,d,k})$ & $(2c,32y_4^2)$ &
$(32y_4^2,\infty)$
\\
$k,\ \sgn A$ & $1,\ +$ & $2,\ -$ & $\sgn A_2$ &
\multicolumn{2}{|c|}{$+$} & $k,\ \sgn A$ & $1,\ +$ & $2,\ -$
\\ \hline
\multicolumn{3}{|c|}{$0<y_3<x_3$} & \multicolumn{3}{|c|}{} &
\multicolumn{3}{|c|}{$0<x_1<y_1$}
\\ \cline{1-9}
\multicolumn{3}{|c|}{$(c,d)\in \ell_3$} &
\multicolumn{3}{|c|}{$(c,d)\in \Omega_4$} &
\multicolumn{3}{|c|}{$(c,d)\in \ell_4$}
\\
\multicolumn{3}{|c|}{Две полубесконечные} &
\multicolumn{3}{|c|}{Семейство грушевидных} &
\multicolumn{3}{|c|}{Грушевидная}
\\
\multicolumn{3}{|c|}{поверхности Бертрана} &
\multicolumn{3}{|c|}{поверхностей Бертрана} &
\multicolumn{3}{|c|}{поверхность Бертрана}
\\ \hline
$I_{c,d,k}$ & $(-\infty,-y_3)$ & $(-y_3,0)$ & $I_{c,d,k}$ &
$(-\infty,-\sqrt[4]{-d})$ & $(-\sqrt[4]{-d},0)$ & $I_{c,d,k}$ & $(-\infty,-y_4)$ & $(-y_4,0)$
\\
$r_{c,d}(I_{c,d,k})$ & $(0, \infty)$ & $(-\infty,0)$ &
$r_{c,d}(I_{c,d,k})$ & $(-\frac{\pi}{4y_4}, -r_4)$ & $(-r_4,0)$ &
$r_{c,d}(I_{c,d,k})$ & $(\frac{-\pi}{4y_4},-r_{4,1})$ &
$(-r_{4,1},0)$
\\
$F_{c,d}(I_{c,d,k})$ & $(0, \infty)$ & $(\infty,0)$ &
$F_{c,d}(I_{c,d,k})$ & $(0, \frac{1}{2y_4})$ & $(\frac{1}{2y_4},0)$ &
$F_{c,d}(I_{c,d,k})$ & $(0, \frac{1}{2y_4})$ & $(\frac{1}{2y_4},0)$
\\
$R_{c,d}(I_{c,d,k})$ & $(2c,0)$ & $(0,\infty)$ & $R_{c,d}(I_{c,d,k})$
& $(2c,32y_4^2)$ & $(32y_4^2,\infty)$ & $R_{c,d}(I_{c,d,k})$ &
$(2c,16c)$ & $(16c,\infty)$
\\
$k,\ \sgn A$ & $1,\ +$ & $2,\ -$ & $k,\ \sgn A$ & $1,\ +$ & $2,\ -$ &
$k,\ \sgn A$ & $1,\ +$ & $2,\ -$
\\ \hline
\multicolumn{3}{|c|}{$x_3=y_3=x_4=\sqrt{-c/2},$ $y_4=0$} &
\multicolumn{3}{|c|}{$x_4>0,$ $y_4>0$} &
\multicolumn{3}{|c|}{$x_4=0,$ $y_4=y_1=x_1=\sqrt{c/2}$}
\\ \hline
\end{tabular}}
\caption{\label{tab:inter} Интервалы $I_{c,d,k}$ и их образы (см.\
комментарии \ref {com:tab}, \ref {com:2}, \ref {com:1}).}
\end{center}
\end{table}

\begin{comm} \label {com:tab}
В таблице \ref {tab:inter} указаны максимальные интервалы
$I_{c,d,k},$ в которых принимает значения координата $\theta$ на
поверхностях Бертрана (\ref {scdk}), и области значений монотонных
функций $r_{c,d}(\theta),$
$F_{c,d}(\theta):=f_{c,d,k}(r_{c,d}(\theta)),$ $R_{c,d}(\theta)$
(натуральный параметр на меридиане, радиус параллели при $\mu=1,$
скалярная кривизна). В таблице использованы обозначения
$r_{4,1}:=(\pi/2-1)/(4y_4)$ и
$$
r_1:=\frac{\pi}{2}\frac{y_1}{y_1^2-x_1^2}-\frac{\arctg\sqrt{y_1/x_1}}{y_1-x_1},
\qquad
r_4:=\frac{\pi}{8y_4}-\frac{1}{4x_4}\ln\frac{\sqrt{x^2_4+y^2_4}-x_4}{y_4}.
$$
\end{comm}

\begin{corollary} [(Геометрия поверхностей Бертрана)] \label {cor:geom}
{\rm(A)} На любой поверхности Бертрана
$(S,ds^2)=(S_{c,d,k},ds^2_{\mu,c,d})$ {\rm(см.\ (\ref{scdk}))}
существует единственная граничная окружность $\{\hat{r}\}\times S^1,$
определенная условием $\lim_{r\to\hat{r}}f(r)=0$ (т.е.\ стягивающаяся
в точку, называемую {\em полюсом} поверхности), причем
$\lim_{r\to\hat{r}}\theta(r)=-\infty$ в случае $k=1,$ и
$\lim_{r\to\hat{r}}\theta(r)=0$ в случае $k=2.$ В полюсе поверхность
имеет коническую особенность с полным углом
$2\pi\lim_{r\to\hat{r}}|f'(r)|,$ равным $2\pi/\mu$ при $k=1$ или
$\infty$ при $k=2.$ Полюс  является притягивающим центром для любого
замыкающего центрального потенциала.

{\rm(B)} Поверхность $(S,ds^2)=((a,b)\times
S^1,ds^2)\subset(S_{c,d,k},ds^2_{\mu,c,d})$ реализуема в $\RR^3$ как
поверхность вращения тогда и только тогда, когда $|f'(r)|\le1$ всюду
на $S,$ где
$f'(r)=-(\theta+d\theta^{-3})/(\mu\sqrt{\theta^2+c-d\theta^{-2}}).$ В
частности, дополнительная поверхность Бертрана {\rm($k=2$)} всегда
нереализуема; при $d\le0\le c$ основная поверхность Бертрана
{\rm($k=1$)} реализуема тогда и только тогда, когда $\mu\ge1$
(включая стандартную полусферу при $(\mu,c,d)=(1,1,0)$ и евклидову
плоскость при $(\mu,c,d)=(1,0,0)$); при
$(c,d)\in\Omega_2\cup\overline{\Omega_3}\setminus\{(0,0)\}$ основная
поверхность Бертрана {\rm($k=1$)} нереализуема ни при каком $\mu$
(включая стандартную плоскость Лобачевского при
$(\mu,c,d)=(1,-1,0)$).

{\rm(C)} Скалярная кривизна Римана $R=-2f''(r)/f(r)=:R_{c,d}(\theta)$
на поверхности Бертрана {\rm(\ref {scdk})} удовлетворяет соотношениям
 $$
\frac{R_{c,d}(\theta)}2=c-\frac{6d}{\theta^{2}}-\frac{3cd}{\theta^{4}}+\frac{2d^2}{\theta^{6}},
\qquad
\frac{R_{c,d}'(\theta)}{4!}=\frac{d}{\theta^5}\left(\theta^2+c-\frac{d}{\theta^{2}}\right)
= \frac{d}{\mu^2f^2\theta^5},
 $$
т.е.\ $R=R(\theta)$ постоянна при $d=0,$ возрастает при $d<0,$
убывает при $d>0$ и имеет 
следующие предельные значения: в основном полюсе $R\to2c,$ 
в дополнительном полюсе $R\to+\infty,$ на экваторе $R\to32y_4^2(c,d)>0$ при $d\ne0,$  
на основном абсолюте $R\to-2\Delta(c,d)/x_3^{2}(c,d)\le0$ при $(c,d)\ne(0,0),$
на дополнительном абсолюте $R\to-2\Delta(c,d)/y_3^{2}(c,d)\le0.$ 
В частности, $R>0$ при $d\le0\le c$ и
$(c,d)\ne(0,0),$ а также при $(c,d)\in\overline{\Omega_4}$ и $k=2$;
$R<0$ при $c\le0\le d$ и $(c,d)\ne(0,0),$ а также при
$(c,d)\in\overline{\Omega_3}$ и $k=1$; $R$ имеет непостоянный знак
при остальных 
 $(c,d,k),$ таких что $(c,d)\ne(0,0).$

{\rm(D)} На любой поверхности вращения или, более общо, на ее $N-$мерном аналоге, т.е.\ на римановом многообразии $(a,b)\times S^{N-1}$ с римановой метрикой $ds^2=dr^2+f^2(r)d\varphi^2,$
где $N\in\NN$ и через $d\varphi^2$ обозначена стандартная риманова метрика на
$(N-1)-$мерной единичной сфере $S^{N-1},$ оператор Лапласа-Бель\-тра\-ми действует на центральные функции $h=h(r)$ по формуле
 $$
 \Delta h(r)=h''(r)+(N-1)h'(r)f'(r)/f(r)=(f^{N-1}(r)h'(r))'/f^{N-1}(r).
 $$
При $N=3$ обобщенный гравитационно-кулоновский потенциал
$\Theta=\Theta(r),$ определенный условием $\Theta'(r)=1/f^2(r),$
является гармонической функцией.
\end{corollary}

\begin{comm} \label {com:2} Опишем бифуркации поверхностей Бертрана (\ref {scdk})
при движении пары параметров $(c,d)\in\RR^2$ вокруг начала координат.

(A) При $(c,d)\in\Omega_1\cup \ell_4\cup\Omega_4$ поверхность
Бертрана состоит из двух связных компонент, на которых
$\theta\in(-\infty,-(-d)^{1/4})$ и $\theta\in(-(-d)^{1/4},0)$
соответственно. Эти компоненты являются двумя половинками
``грушевидной'' аналитической римановой поверхности
$(-\infty,0)\times S^1,$ разрезанной по единственному экватору (где
$f'=0$). {\it Грушевидность} поверхности вращения означает, что
существует единственный экватор, функции $r$ и $f$ монотонны на
каждой половинке поверхности вне экватора, принимают значения в
конечных интервалах (т.е.\ поверхность ограничена), причем в каждой
половинке есть полюс (т.е.\ $\inf f=0,$ откуда $\sup f$ достигается
на экваторе), и в полюсах поверхность имеет конические особенности с
разными полными углами (равными $2\pi/\mu$ и $\infty$ в основном и
дополнительном полюсах соответственно). При $\mu\in\QQ$ соответствующий замыкающий
(осцилляторный) потенциал на разных половинках пропорционален одной
и той же аналитической функции, однако знаки коэффициентов
пропорциональности различны на разных половинках. Для
каждого знака потенциала все неособые ограниченные орбиты замкнуты
(т.е.\ на всей поверхности без полюсов потенциал является замыкающим,
полулокально и локально замыкающим, но не слабо и не сильно
замыкающим). Длина каждой замкнутой геодезической, образованной двумя
меридианами, равна $\mu L/2,$ а при $\mu\in\QQ$ все
остальные геодезические замкнуты и проходят через экватор, причем все
геодезические кроме экватора и меридианов имеют одну и ту же длину $qL,$ где $L$ --- длина экватора, $\mu/2=:q/p$ --- несократимая дробь. Эти грушевидные поверхности 
являются примерами {\it поверхностей Таннери}~\cite[теорема~4.13 и приложение~А]{Besse}, т.е.\ поверхностей вращения, все геодезические которых замкнуты. При $\mu=2$ получаем примеры {\it поверхностей Цолля}~\cite[следствие~4.16]{Besse}, т.е.\ поверхностей, все геодезические которых имеют одну и ту же длину.

При приближении точки $(c,d)\in\Omega_1$ к лучу $\ell_1$
дополнительная половинка ``исчезает'' (вырождается в экватор), а
основная превращается в проколотую полусферу (при $\mu=1$) или в
локально изометричную ей поверхность (при произвольном $\mu>0$).
Максимальное аналитическое продолжение полусферы --- сфера --- сходна
с грушевидными поверхностями, с одним лишь отличием --- она
симметрична относительно экватора. Так как осцилляторный потенциал
равен $\theta^{-2},$ он имеет особенность на экваторе сферы (т.е.\
при $\theta=0$), а потому не продолжим на всю сферу. Гравитационный
потенциал на ``основной'' полусфере равен $\theta,$ а его
аналитическое продолжение на сферу (как и для осцилляторного
потенциала на грушевидных поверхностях) является замыкающим,
полулокально и локально замыкающим, но не слабо и не сильно
замыкающим.

При приближении точки $(c,d)\in\Omega_4$ к дуге $\ell_3$ экватор
бесконечно удлиняется и отдаляется от полюсов (``превращается в
абсолют''), в результате чего половинки ``отделяются друг от друга'',
и грушевидная поверхность распадается на две полубесконечные
поверхности (см.\ (C) ниже).

(B) При $(c,d)\in\Omega_2$ поверхность связна (т.е.\ состоит из одной
лишь половинки --- основной) и \textit{полубесконечна}, т.е.\ $r$ и
$f$ возрастают (как функции друг от друга), причем имеются полюс
($\inf r>-\infty,$ $\inf f=0$) и ``абсолют'' ($\sup r=\sup
f=+\infty$). Тем самым, каждая точка находится на конечном расстоянии
от полюса и может быть удалена от него на любое расстояние. В полюсе
поверхность имеет коническую особенность с полным углом $2\pi/\mu.$
Все геодезические с $K\ne0$ незамкнуты, имеют бесконечную длину (в
обе стороны).

При приближении точки $(c,d)\in\Omega_2$ к лучу $\ell_1$ или $\ell_2$ поверхность превращается 
в проколотую полусферу или плоскость Лобачевского (в случае $\mu=1$). 

(C) При $(c,d)\in\Omega_3\cup \ell_3$ поверхность состоит из двух
связных компонент --- основной и дополнительной, причем каждая
компонента полубесконечна (см.\ выше), и компоненты друг другу не
изометричны (и не подобны), даже локально. В полюсах поверхность
имеет конические особенности с разными полными углами (равными
$2\pi/\mu$ и $\infty$ в основном и дополнительном полюсах
соответственно). Все геодезические с $K\ne0$ имеют бесконечную длину.

При приближении точки $(c,d)\in\Omega_3$ к лучу $\ell_2$
дополнительная компонента ``исчезает'', а основная превращается в
проколотую плоскость Лобачевского (при $\mu=1$) или в локально
изометричную ей поверхность (при любом $\mu>0$).
\end{comm}


\begin{corollary} [(Классификация поверхностей Бертрана)] \label {cor:class}
Рассмотрим четверки чисел $(\mu,c,d,\theta_0),$ такие что
$\mu,c,d\in\RR,$ $\mu>0,$ $\theta_0\in I_{c,d},$ и пары чисел
$(\lambda,\phi)\in\RR_{>0}\times[0,2\pi],$ такие что $\lambda=1$ в
случае $\gamma(\phi)=(0,0).$ На поверхности Бертрана
$(S_{c,d},ds^2_{\mu,c,d})$ рассмотрим координаты
$\theta,\varphi\mod2\pi$ {\rm(см.\ (\ref {scdk}))}, меридианы
$I_{c,d}\times\{\varphi_0\}\subset S_{c,d}$ и параллели
$\{\theta_0\}\times S^1\subset S_{c,d}.$

{\rm(A)} Для различных четверок $(\mu_i,c_i,d_i,\theta_i),$ $i=1,2,$
поверхности Бертрана $(S_{c_i,d_i},ds^2_{\mu_i,c_i,d_i})$ не
изометричны друг другу (с сохранением семейства меридианов) ни в
каких окрестностях параллелей $\{\theta_i\}\times S^1\subset
S_{c_i,d_i}.$

Поверхность $(S_{c,d},ds^2_{\mu,c,d})$ локально изометрична
поверхности $(S_{c,d},ds^2_{1,c,d}).$ Точнее, окрестность любого
меридиана поверхности $(S_{c,d},ds^2_{\mu,c,d})$ изометрично
вкладывается в поверхность $(S_{c,d},ds^2_{1,c,d})$ посредством
сопоставления $(\theta,\varphi)$ $\mapsto(\theta,\varphi/\mu).$ Для
различных троек $(c_i,d_i,\theta_i),$ $i=1,2,$ поверхности
$(S_{c_i,d_i},ds^2_{1,c_i,d_i})$ не являются локально изометричными
друг другу (с сохранением семейства меридианов) ни в каких
окрестностях точек $(\theta_i,0)\in S_{c_i,d_i}.$

{\rm(B)} Поверхность $(S_{c,d},ds^2_{\mu,c,d})$ подобна поверхности
$(S_{\gamma(\phi)},ds^2_{\mu,\gamma(\phi)})$ и локально подобна
поверхности $(S_{\gamma(\phi)},ds^2_{1,\gamma(\phi)})$ с
коэффициентом подобия $1/\lambda,$ где пара $(\lambda,\phi)$
однозначно определяется соотношениями $c=\lambda^2c_\phi,$
$d=\lambda^4d_\phi.$ Точнее, поверхность
$(S_{c,d},\lambda^2ds^2_{\mu,c,d}),$ подобная поверхности
$(S_{c,d},ds^2_{\mu,c,d}),$ изометрично отображается на
$(S_{\gamma(\phi)},ds^2_{\mu,\gamma(\phi)})$ сопоставлением
$(\theta,\varphi)\mapsto(\theta/\lambda,\varphi).$ Для различных
троек $(\mu_i,\phi_i,\theta_i),$ $i=1,2,$ поверхности
$(S_{\gamma(\phi_i)},ds^2_{\mu_i,\gamma(\phi_i)})$ не подобны друг
другу (с сохранением семейства меридианов) ни в каких окрестностях
параллелей $\{\theta_i\}\times S^1\subset S_{\gamma(\phi_i)}.$ Для
различных пар $(\phi_i,\theta_i),$ $i=1,2,$ поверхности
$(S_{\gamma(\phi_i)},ds^2_{1,\gamma(\phi_i)})$ не являются локально
подобными друг другу (с сохранением семейства меридианов) ни в каких
окрестностях точек $(\theta_i,0)\in S_{\gamma(\phi_i)}.$

{\rm(C)} Поверхность ${(S_{c,d},ds^2_{\mu,c,d})}$ проективно
диффеоморфна области на одной из трех поверхностей
${(S_{c_0,d_0},ds^2_{\mu,c_0,d_0})}$ для $d_0:=\sgn
d\in\{-1,0,1\}$ и любого $c_0\ge c_0(c,d),$ где $c_0(c,d):=0$
при $d=0,$ $c_0(c,d):=-2$ при $d<0,$
$c_0(c,d):=c/\sqrt{d}$ при $d>0.$ Точнее, имеется
вложение ${(S_{c,d},ds^2_{\mu,c,d})}\hookrightarrow{(S_{c_0,d_0},ds^2_{\mu,c_0,d_0})},$
переводящее (непараметризованные) геодезические в
гео\-де\-зи\-чес\-кие и задаваемое следующим сопоставлением. При
$d=0$ оно задается сопоставлением
$(\theta,\varphi)\mapsto(\theta,\varphi),$ причем при $c\ge0$ оно
биективно, а при $c<0$ переводит поверхность (являющуюся при $\mu=1$ проколотой плоскостью Лобачевского) на проколотый круг
$(-\infty,-\sqrt{-c})\times S^1\subset S_{c_0,0}.$ При $d\ne0$ оно
задается сопоставлением
$(\theta,\varphi)\mapsto(\theta/|d|^{1/4},\varphi),$ причем
при $(c,d)\in\Omega_1\cup\ell_4\cup\Omega_4\cup\ell_3$ оно биективно, при
$(c,d)\in\Omega_3$ отображает поверхность на объединение двух проколотых кругов
$(-\infty,-x_3(c,d)/(-d)^{1/4})\times S^1$ и $(-y_3(c,d)/(-d)^{1/4},0)\times
S^1,$ а при $d>0$ отображает поверхность на проколотый круг
$(-\infty,-x_3(c,d)/d^{1/4})\times S^1.$
\end{corollary}

\begin{comm}\label {com:1}
Поверхности Бертрана (\ref {scdk}) являются аналитическими, попарно
неизометричными (даже локально, с сохранением действия группы
вращений) поверхностями и образуют гладкое трехпараметрическое
семейство с параметрами $(\mu,c,d)$ или $(\mu,\lambda,\phi).$ Два
параметра $\mu,\lambda>0$ имеют простой геометрический смысл --- при
их изменении (с подходящим линейным масштабированием координатных осей) компоненты римановой метрики умно\-жа\-ют\-ся на константы, а
потому полученная поверхность локально подобна 
исходной, с сохранением семейства мери\-ди\-а\-нов. При классификации метрик (\ref {scdk})
с точностью до изометрии (соответственно локальной изометрии,
подобия, локаль\-но\-го подобия), сохраняющей семейство
мери\-ди\-а\-нов, полным инвариантом является набор
параметров $(\mu,c,d)$ (соответственно $(c,d),$ $(\mu,\phi),$
$\phi$), а остальные параметры несущественны для классификации.
Класс проективной (соот\-вет\-ст\-венно локальной проективной)
диффеоморфности поверхности (\ref {scdk})
вполне определяется параметром $\mu>0$ и знаком $\sgn d$
(соот\-вет\-ст\-венно знаком $\sgn d$). 
Свойства проективно эквивалентных метрик исследуются, например, в~\cite{Matv}.
\end{comm}

\subsection {Уравнения орбит движения точки по поверхности Бертрана в поле замыкающего потенциала} \label {subsec:orbits}

Выпишем в явном виде уравнения неособых орбит движения точки по
поверхностям Бертрана (\ref {scdk}) (т.е.\ (\ref {f}) и (\ref {eq:sol})).

Поверхность первого типа, $f(r)=\xi f_{c}(r),$ неособая орбита
$\theta=\theta_{c,0}(r(\varphi))$ для гравитационного потенциала
$V(r)=V_{c,0,1}(r)=-A|\theta|+B$ при $A>0$:
$$
 \theta = \pm\frac{\xi^2A}{K^2}\left(1+\sqrt{1+2\frac{E_1K^2}{\xi^2A^2}}\sin(\xi(\varphi-\varphi_0)) \right),
 \qquad \Phi=2\pi/\xi.
$$
Поверхность первого типа, $f(r)=\xi f_{c}(r),$ неособая орбита
$\theta=\theta_{c,0}(r(\varphi))$ для осцилляторного потенциала
$V(r)=V_{c,0,2}(r)=\frac{A}{2\theta^2}+B$ при $A>0$:
$$
 \theta^2 = \frac{\xi^2E_1}{K^2} \left( 1+\sqrt{1-\frac{AK^2}{\xi^2E_1^2}} \sin(2\xi(\varphi-\varphi_0)) \right),
 \qquad \Phi=\pi/\xi.
$$
Поверхность второго типа, $f(r)=f_{c,d}(r)/\mu,$ неособая орбита
$\theta=\theta_{c,d}(r(\varphi))$ для осцилляторного потенциала
$V(r)=V_{c,d,2}(r)=\frac{A}{2\theta^2}+B$ при $(\theta^4+d)A>0$:
$$
 \theta^2 = \frac{E_1}{\mu^2K^2} \left( 1+\sqrt{1-\frac{\mu^2K^2}{E_1^2} \left( A-\mu^2K^2d\right)}
 \sin\left(2\frac{\varphi-\varphi_0}\mu\right)\right),
 \qquad \Phi=\pi\mu.
$$
Здесь $E=({\dot
r}^2+f^2(r){\dot\varphi}^2)/2+V(r)=K^2(\mu^4(d\theta/d\varphi)^2+1/f^2(r))/2+V(r)$
--- константа полной энергии, $K\ne0$ --- интеграл кинетического момента,  $E_1=E-B-\mu^2K^2c/2,$ $c,d\in\RR,$ $\xi=1/\mu,$ $\mu$ --- положительное раци\-о\-наль\-ное число, определяемое метрикой (см.\ (\ref {metric}), (\ref {f}), (\ref {eq:sol})),
$\Phi$ --- минимальный поло\-жи\-тель\-ный период периодических
функций $r(\varphi),$ отличных от констант.

\subsection {Обобщение технической теоремы \ref {bert2} Бертрана}

{\it Обобщенным семейством уравнений Бертрана} назовем
однопараметрическое семейство дифференциальных уравнений
$\frac{d^{2}z}{d\varphi^{2}}+\rho(z)=\frac{1}{K^2}\Psi(z)$ на
интервале $(a,b) \subset \mathbb{R}^1$ с параметром $K \in
\RR\setminus\{0\},$ где $\Psi(z)$ и $\rho(z)$ --- функции класса
$C^{\infty},$ определенные на интервале $(a,b).$ Следующее
определение аналогично определению \ref {closing}.

\begin{definition} \label {def:rho:closing}
Функцию $\Psi=\Psi(z)$ на интервале $(a,b)$ будем называть
\textit{замыкающей для функции $\rho=\rho(z)$} (или {\em
$\rho$-замыкающей}), если
\begin{enumerate}
\item[$(\exists)$] существует значение параметра $K=\hat K\in\RR\setminus\{0\},$
при котором уравнение имеет ограниченное непостоянное решение $\hat z=\hat
z(\varphi)$;
\item[$(\forall)$]
все ограниченные непостоянные решения $z=z(\varphi)$ уравнения со всевозможными
значениями параметра $K$ являются периодическими функциями с попарно
соизмеримыми периодами.
\end{enumerate}
Функцию $\Psi(z)$ будем называть \textit{локально замыкающей для
функции $\rho(z)$} (или {\em локально $\rho$-замыкающей}), если
\begin{enumerate}
\item[$(\exists)$\loc] существует значение параметра $K=K_0,$ при котором
уравнение имеет невырожденное устойчивое положение равновесия
$z_0\in(a,b)$;
\item[$(\forall)$\loc] для любой пары $(K_0,z_0)\in(\RR\setminus\{0\})\times(a,b),$ удовлетворяющей
условию $(\exists)$\loc, существуют $\varepsilon,\delta>0,$ такие что
все ограниченные непостоянные решения $z=z(\varphi)$ уравнений со
значениями параметра $K\in(K_0-\delta,K_0+\delta),$ такие что
$z(\RR^1)\subseteq[z_0-\varepsilon,z_0+\varepsilon],$ являются
периодическими функциями с попарно соизмеримыми периодами.
\end{enumerate}
Функцию $\Psi(z)$ будем называть \textit{полулокально замыкающей для
функции $\rho(z)$} (или {\em полулокально $\rho$-замыкающей}), если
выполнены условия $(\exists),$ $(\forall)$\loc\ и следующее условие:
\begin{enumerate}
\item[$(\forall)$\sloc] все ограниченные непостоянные решения $z=z(\varphi)$ уравнения
при $K=\hat K,$ такие что $z(\RR^1)\subseteq\hat z(\RR^1),$ являются
периодическими функциями с попарно соизмеримыми периодами, где $\hat
K$ и $\hat z=\hat z(\varphi)$ -- значение параметра и решение из
$(\exists).$
\end{enumerate}
Функцию $\Psi(z)$ назовем \textit{сильно} (соответственно
\textit{слабо}) \textit{$\rho$-замыкающей}, если любая точка
$z_0\in(a,b)$ является невырожденным устойчивым (соответственно
устойчивым) положением равновесия уравнения при некотором $K=K_0,$
зависящем от $z_0,$ и выполнено условие $(\forall)$\loc\
(соответственно его аналог для всякой пары
$(K_0,z_0)\in(\RR\setminus\{0\})\times(a,b),$ такой что $z_0$ ---
устойчивое положение равновесия уравнения при $K=K_0$).
\end{definition}

Доказательства теорем \ref {our1}, \ref {our_1_1}, \ref {statement}
основаны на следующем обобщении технической теоремы \ref{bert2}
Бертрана, путем замен $z(r)-\zeta=-\Theta(r)=-\mu^2\theta(r),$ где
$\Theta'(r)=1/f^2(r),$ $\rho(z(r))=f'(r)/f(r),$
$\Psi(z)=f^2(r)V'(r)=-dV(r(z))/dz,$ $D=\mu^8d.$

\begin{theorem} [(Обобщенная техническая теорема Бертрана)]\label{genbert}
Рассмотрим однопараметрическое семейство дифференциальных уравнений
$\frac{d^{2}z}{d\varphi^{2}}+\rho(z)=\frac{1}{K^2}\Psi(z)$ на
интервале $(a,b) \subset \mathbb{R}^1$ с параметром $K \in
\RR\setminus \{0\},$ где $\Psi=\Psi(z)$ и $\rho=\rho(z)$ --- функции
класса $C^{\infty},$ определенные на интервале $(a,b).$ Если $\Psi$
является  полулокально $\rho$-замыкающей (или $\rho$-замыкающей или
сильно или слабо $\rho$-замыкающей), то она является локально
$\rho$-замыкающей.

Пусть функция $\rho$ не имеет нулей на интервале $(a,b).$ Тогда в
интервале $(a,b)$ классы замыкающих, полулокально замыкающих,
локально замыкающих, сильно замыкающих и слабо замыкающих для $\rho$
функций $\Psi$ совпадают, причем существует не более двух замыкающих
функций $\Psi(z)$ с точностью до положительной мультипликативной константы и эти функции определяются следующими условиями:

{\rm(a)} если $\rho'|_{(a,b)} = \const>0,$ то существуют ровно две (с
точностью до положительной мультипликативной константы)
$\rho-$замыкающие функции $\Psi$ на $(a,b),$ а именно
$\Psi_i(z)=A_i/\rho^{i^2-1}(z),$ $i=1,2$
(т.е.\ функции, отвечающие обобщенному гравитационному и
осцилляторному законам сил на соответствующей поверхности вращения),
где $A_i\ne0$ --- произвольная мультипликативная константа, такая что
$A_i\rho^i(z)>0,$
причем минимальный положительный период любого ограниченного непостоянного решения
равен $\Phi_i=2\pi/(i\sqrt{\rho'})$;

{\rm(b)} если $\rho|_{(a,b)}$ является рациональной функцией вида
$\rho(z) = \frac{(z-\zeta)^4+D}{\mu^2(z-\zeta)^3},$ где $D=\const
\neq 0,$ $\mu=\const>0,$ $\zeta=\const\notin(a,b),$ то существует
единственная (с точностью до положительной мультипликативной
константы) $\rho-$замыкающая функция на $(a,b):$ $
\Psi(z)=\Psi_2(z)=\frac{A}{(z-\zeta)^3} $ (т.е.\ отвечающая
осцилляторному закону сил на соответствующей поверхности вращения),
где $A\ne0$ --- произвольная мультипликативная константа, такая что
$A((z-\zeta)^4+D)>0,$ причем минимальный положительный период любого
ограниченного непостоянного решения
равен $\Phi = \pi\mu$;

{\rm(c)} если $\rho(z)$ не имеет ни одного из указанных выше видов,
то не существует $\rho-$замыкающих функций на $(a,b).$

В случаях $(a)$ и $(b)$ каждая точка $z \in (a,b)$ является
невырожденным устойчивым положением равновесия уравнения при
$K=K_i:=\pm\sqrt{{A_i}/{\rho^{i^2}(z)}},$ $i=1,2$
(в случае {\rm(a)}) и $K=\pm\mu\sqrt{\frac{A}{(z-\zeta)^4+D}}$ (в
случае {\rm(b)}), а при других значениях параметра $K$ не является
положением равновесия.
\end{theorem}

Отметим, что в теореме \ref{genbert} (в отличие от теоремы
\ref{bert2}) не требуется аналитичность функций $\Psi_i(z)$ и
$\rho(z),$ а постоянная $\mu>0$ не обязана быть рациональной
(поскольку в ней не требуется, чтобы все периоды были соизмеримы с
$2\pi,$ а требуется лишь попарная соизмеримость периодов, как в
замечании \ref {rem:xi}).

\section{Частный случай: конус} \label {sec:cone}

Согласно следствию \ref{cone}, теорема \ref{bert1} Бертрана обобщается на семейство ``рациональных'' конусов. Прочие конусы (``иррациональные'') не допускают обобщение
теоремы Бертрана.
Поясним, какую роль играет условие рациональности константы $\xi$ в
условии следствия \ref{cone}. Рассмотрим конструкцию, описанную после
теоремы \ref{our1}: разрежем конус $S$ по образующей и развернем.
Получим некоторый сектор на евклидовой плоскости. Заметим, что угол
при вершине (равный $2\pi\xi$) может быть и больше $2\pi,$ в этом
случае конус не вложится в $\mathbb{R}^{3}$ как поверхность вращения,
но следствие \ref{cone} останется верным. Далее рассмотрим
поверхность $\tilde{S},$ являющуюся одновременно разветвленным
накрытием конуса $S$ и разветвленным накрытием евклидовой плоскости,
где количество листов любого из накрытий минимально возможное.
Поверхность $\tilde{S}$ можно построить следующим образом. Будем
накладывать на плоскость сектора, полученные разворотом разрезанного
по образующей конуса, каждый следующий поворачивая так, чтобы берег
каждого следующего сектора совпал с противоположным берегом
предыдущего. Таким образом получим разветвленное накрытие $\tilde{S} \to \RR^2$ над
плоскостью. При этом, если конус был
``рационален'', т.е. угол при вершине был соизмерим с $2\pi$ (иными
словами, $\xi \in \QQ,$ обозначим $\xi=\frac{p}{q}),$ то через $q$
шагов построения $\tilde{S}$ берег очередного сектора совпадет с
берегом первого сектора; в этом случае прекратим процесс построения
$\tilde{S}$ и накрытие $\tilde{S} \to \RR^2$ будет $p-$листным. В
противном случае берега разных секторов никогда не совпадут и
накрытие будет бесконечнолистным. Замощенная плоскость и образ
траектории движения точки по поверхности $\tilde{S}$ изображены на рис.~\ref{ris:1}.

\begin{figure}[t]
\begin{center}
$\includegraphics[width=90mm]{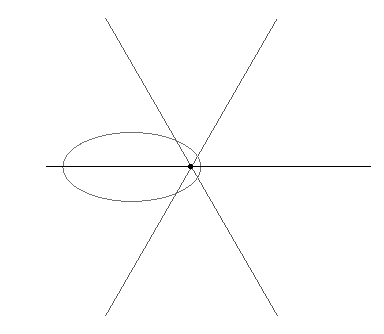}$
\caption{Образ траектории движения по конусу на замощенной плоскости;
$\xi={1}/{6}.$} \label{ris:1}
\end{center}
\end{figure}

Таким образом, поверхность $\tilde{S}$ накрывает исходный конус
$q-$листно, если $\xi$ рационально и равно $\frac{p}{q},$
бесконечнолистно в противном случае. Накрытие строится естественным
образом: его листами будут сектора, из которых составлена поверхность
$\tilde{S}$; каждый такой сектор --- развертка конуса, поэтому можно
определить отображение $\tilde{S} \to S,$ переводящее в точку с
координатами $(r, \varphi)$ на конусе все точки поверхности
$\tilde{S},$ имеющие те же координаты на том секторе, которому
принадлежат. Наглядно это может быть представлено следующим образом.
Возьмем сектора, полученные разворачиванием разрезанного по
образующей конуса, и расположим их над разверткой конуса ``друг над
другом'' --- $q$ экземпляров, если $\xi$ рационально, счетное
множество в противном случае; противоположные края лежащих друг над
другом секторов считаются склеенными. Точки этих секторов --- это
точки поверхности $\tilde{S},$ отображение $\tilde{S} \to S$ задается
естественной проекцией. Образ траектории движения по конусу на
поверхности $\tilde{S},$ представленной как пояснено выше,
продемонстрирован на рис. \ref{ris:2}. На рис. \ref{ris:3} приведен
пример замкнутой траектории на ``рациональном'' конусе, вид сверху.

\begin{figure}[t]
\begin{center}
$\includegraphics[width=90mm]{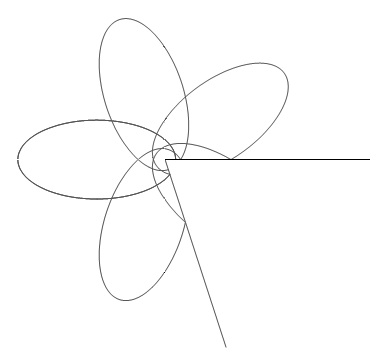}$
\caption{Образ траектории движения по конусу на поверхности $\tilde{S},$ представленной объединением секторов, расположенных ``друг над другом'' на плоскости (на каждом секторе своя дуга траектории); $\xi={5}/{6}.$}
\label{ris:2}
\end{center}
\end{figure}

\begin{figure}[t]
\begin{center}
$\includegraphics[width=110mm]{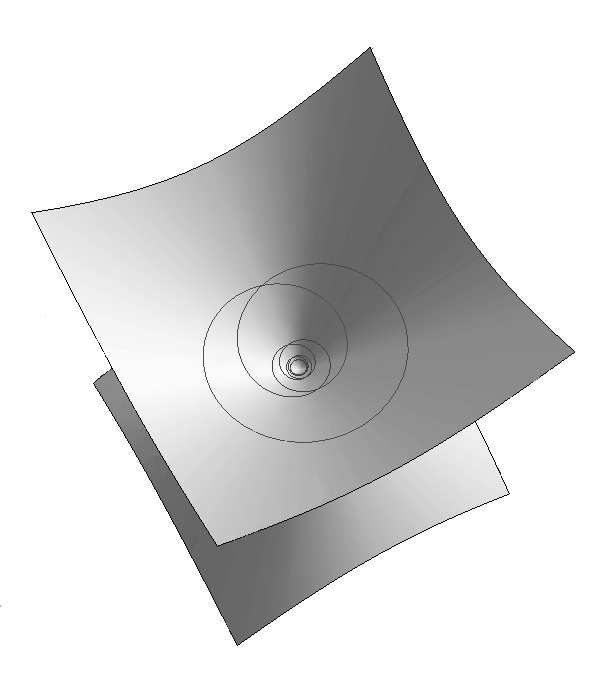}$
\caption{Замкнутая траектория на рациональном конусе;
$\xi=1/8.$} \label{ris:3}
\end{center}
\end{figure}

\section{Доказательство обобщенной технической теоремы \ref {genbert} Бертрана} \label {sec:technical}

Обобщенное уравнение Бертрана имеет вид
$z''_{\varphi\varphi}=-U'(z),$ где $U(z):=\int(\rho(z)-\Psi(z)/K^{2})dz.$ Поэтому его решения удовлетворяют закону сохранения
энергии $(z')^2/2+U(z)=:E,$ а потому оно может быть сведено к уравнению первого порядка
$(z')^2=2E-2U,$ где $E$ --- константа (называемая уровнем энергии в
обобщенной задаче Бертрана), зависящая от решения, т.е.\ к уравнению
$$
(z')^2=R(z)+\C W(z)+2E=2E-2U_\C(z),
$$
где $R(z):=-2\int\rho(z)dz,$ $\C:=2/{K^2},$ $W(z) :=
\int\Psi(z)dz,$ $U_\C(z):=-(R(z)+\C W(z))/2.$ Величина $K^2E$ совпадает с
{\em уровнем энергии} (см.\ \S\ref{subsec:orbits}), а функция
$K^2U_\C(z)$ --- с {\em эффективным потенциалом} (см.\ определение
\ref {traj}(a)) в задаче о движении точки по поверхности $S$ в
центральном поле с потенциалом $V=V(r).$

Прежде всего докажем три вспомогательных утверждения о свойствах
решений дифференциального уравнения $z''_{\varphi\varphi}=-U'(z),$
где $U=U(z)$ \--- гладкая функция на интервале $(a,b).$

\begin{propos} \label {pro:1}
Если функция $z=z(\varphi)$ является решением уравнения
$z''_{\varphi\varphi}=-U'(z)$ и $z'(\varphi_0)=0,$ то график функции
$z=z(\varphi)$ симметричен относительно прямой
$\{\varphi=\varphi_0\},$ иными словами
$z(\varphi)=\tilde{z}(\varphi):=z(\varphi_0-(\varphi-\varphi_0)).$
\end{propos}

\begin{proof}
То, что функция $\tilde{z}(\varphi)$ удовлетворяет уравнению
$z''_{\varphi\varphi}=-U'(z),$ доказывается прямой подстановкой. А
поскольку $z(\varphi_0)=\tilde{z}(\varphi_0)$ и
$z'(\varphi_0)=\tilde{z}'(\varphi_0),$ то решения $z(\varphi)$ и
$\tilde{z}(\varphi)$ совпадают. Утверждение доказано.
\end{proof}

\begin{propos} \label {pro:2.5}
Пусть $a<a'<b'<b$ и $E'\in\RR.$ Тогда следующие условия равносильны:

(a) существует ограниченное решение $z(\varphi)$ уравнения
$z''_{\varphi\varphi}=-U'(z)$ с уровнем энергии $E,$ такое что
$a'=\inf z(\RR^1),$ $b'=\sup z(\RR^1)$;

(b) $U(a')=U(b')=E'$ и $U|_{(a',b')}<E'.$

Если $\hat z\in\{a',b'\}$ и выполнено условие (a), то соотношения
$U'(\hat z)\ne0$ и $\hat z\in z(\RR^1)$ равносильны.
\end{propos}

\begin{propos} \label {pro:3}
Пусть $a<a'<b'<b,$ $U(a')=U(b')=E',$ $U|_{(a',b')}<E',$ $E_0:=\min
U|_{[a',b']}.$ Тогда следующие условия равносильны:

(a) для любого $E\in(E_0,E']$ любое ограниченное решение
$z_E(\varphi)$ уравнения $z''_{\varphi\varphi}=-U'(z)$ с уровнем
энергии $E,$ такое что $z_E(0)\in(a',b'),$ периодично;

(b) существует отрезок $[c_1,c_2]\subset (a',b'),$ такой что
$U'|_{[a',c_1)}<0,$ $U'|_{[c_1,c_2]}=0$ и $U'|_{(c_2,b']}>0.$

При выполнении этих условий минимальный положительный период решения
$z_E(\varphi)$ равен
 \begin{equation} \label {eq:3.5_0}
\Phi(E)=2\int_{z_1(E)}^{z_2(E)}\frac{dz}{\sqrt{2E-2U(z)}},
 \end{equation}
где значения $z_1=z_1(E)\in[a',c_1)$ и $z_2=z_2(E)\in(c_2,b']$
определены условиями $U(z_1)=U(z_2)=E.$ Функция $\Phi=\Phi(E)$
непрерывна на полуинтервале $(E_0,E'].$ Если $U''(c_1)=0,$ то
 \begin{equation} \label {eq:3.5_1}
\lim_{E\to E_0}\Phi(E)=\infty.
 \end{equation}
\end{propos}

\begin{proof}[предложений \ref {pro:2.5} и \ref {pro:3}] {\it Шаг 1.}
Предположим, что выполнено какое-либо из условий (a) и (b)
предложения \ref {pro:2.5}.
Тогда существует точка $z_0\in(a',b'),$ такая что $U(z_0)<E'.$
Фиксируем любое число $E\in(U(z_0),E'].$ Пусть $z(\varphi)$ --
локальное решение дифференциального уравнения
$z''_{\varphi\varphi}=-U'(z),$ такое что $z(0)=z_0,$
$z'(0)=\sqrt{2E-2U(z_0)}.$ Тогда уровень энергии на решении
$z(\varphi)$ равен $E.$ Пусть $(z_1,z_2) \subseteq (a,b)$ \---
максимальный интервал по включению, содержащий точку $z_0,$ на
котором $U(z)<E.$ Тогда решение $z(\varphi)$ может быть продолжено на
интервал $(\varphi_1,\varphi_2)$ (с помощью закона сохранения
энергии), где
$\varphi_i:=\int_{z_0}^{z_i}\frac{dz}{\sqrt{2E-2U(z)}},$ $i=1,2.$
Отсюда следует, что при выполнении какого-либо из условий (a) и (b)
предложения \ref {pro:2.5} интервал $(a',b')$ содержит максимальный
по включению интервал, содержащий точку $z_0,$ на котором $U(z)<E'.$
Значит, $(z_1,z_2)\subseteq(a',b')$ (в силу $E\le E'$).

Покажем, что решение $z(\varphi)$ ограничено и
$z(\RR^1)\subseteq[z_1,z_2]$ (т.е.\ $z_1=\inf z(\RR^1)$ и $z_2=\sup
z(\RR^1)$). (Отсюда следует равносильность условий (a) и (b)
предложения \ref {pro:2.5}.) Если $\varphi_1=-\infty$ и
$\varphi_2=\infty,$ то требуемое утверждение очевидно, так как
построенное решение определено на всем $\RR^1$ и $z(\RR^1)=(z_1,z_2)$
(а потому решение непериодично). Осталось рассмотреть случай, когда
для некоторого $i\in\{1,2\}$ выполнено $|\varphi_i|<\infty.$ В этом
случае в силу включения $(z_1,z_2)\subseteq[a',b']\subset(a,b)$
решение может быть продолжено в точку $\varphi_i,$ откуда
$z(\varphi_i)=z_i,$ $U(z_i)=E$ и $z'(\varphi_i)=0.$ Если
$\varphi_1=-\infty$ и $\varphi_2<+\infty,$ то в силу предложения \ref
{pro:1} решение ограничено и
$z(\RR^1)=z((-\infty,\varphi_2])=(z_1,z_2]$ (а потому решение
непериодично). Аналогично показывается, что если $\varphi_1>-\infty$
и $\varphi_2=+\infty,$ то решение ограничено и
$z(\RR^1)=z([\varphi_1,\infty))=[z_1,z_2)$ (а потому решение
непериодично), а в случае $\varphi_1>-\infty$ и $\varphi_2<+\infty$
решение периодично и $z(\RR^1)=z([\varphi_1,\varphi_2])=[z_1,z_2],$
причем минимальный положительный период равен
$\Phi=2\int_{z_1}^{z_2}\frac{dz}{\sqrt{2E-2U(z)}}.$ Таким образом,
решение всегда ограничено и $z(\RR^1)\subset[z_1,z_2],$ что и
требовалось.

Докажем, что для любого $i\in\{1,2\}$ соотношения $U'(z_i)\ne0$ и
$z_i\in z(\RR^1)$ равносильны. Пусть для определенности $i=2.$
Предположим, что $U'(z_2)=0,$ т.е.\ точка $z_2$ является положением
равновесия. Если $z_2\in z(\RR^1)$ (т.е.\ $\varphi_2<\infty,$ см.\
выше), то оба решения $z(\varphi)$ и $z_2(\varphi)\equiv z_2$
удовлетворяют одной системе начальных условий $z(\varphi_2)=z_2,$
$z'(\varphi_2)=0,$ а потому совпадают, что противоречит условию
$z_1<z_2.$ Предположим теперь, что $U'(z_2)\ne 0.$ Тогда
$U|_{[z_0,z_2]}(z)\le E-c(z_2-z)$ для некоторого $c>0,$ откуда
$\varphi_2\le(z_2-z_0)\sqrt{2/c}<\infty,$ а потому $z_2\in z(\RR^1)$
(см.\ выше). Предложение \ref {pro:2.5} полностью доказано.

{\it Шаг 2.} Если выполнено условие (b) предложения \ref {pro:3}, то
$U'(z_1)\ne0$ и $U'(z_2)\ne0,$ а значит, по предложению \ref
{pro:2.5} образ $z(\RR^1)$ решения $z(\varphi)$ является отрезком,
откуда решение периодично (см.\ шаг 1). Отсюда и из предложения \ref
{pro:2.5} следует, что условие (b) предложения \ref {pro:3} влечет
условие (a) предложения \ref {pro:3}.

{\it Шаг 3.} Предположим, что выполнено условие (a) предложения \ref
{pro:3}. Так как $a'<b'$ и $U(a')=U(b'),$ то множество
$A:=\{z\in(a',b')\mid U'(z)=0\}$ непусто. Обозначим $E^*:=\sup U|_A,$
$B:=\{z\in[a',b']\mid U(z)\le E^*\},$ $c_1:=\inf B,$ $c_2:=\sup B.$
Тогда
 \begin{equation} \label {eq:good}
 U'|_{[a',c_1)}<0, 
 \quad U'|_{(c_2,b']}>0, 
\quad E_0\le U|_{[c_1,c_2]}\le U(c_1)=U(c_2)=E^*<E'.
 \end{equation}
Предположим, что $U([c_1,c_2])\ne E^*.$ Тогда $c_1<c_2$ и, в силу
(\ref {eq:good}), найдется интервал
$(z_1^*,z_2^*)\subseteq(c_1,c_2),$ такой что
$U|_{(z_1^*,z_2^*)}<E^*,$ $U(z_1^*)=U(z_2^*)=E^*$ и один из концов
$z_i^*$ ($i\in\{1,2\}$) этого интервала принадлежит $A$ (т.е.\
выполнено $U'(z_i^*)=0$). Пусть теперь $z^*(\varphi)$ -- ограниченное
решение с уровнем энергии $E^*,$ отвечающее интервалу $(z_1^*,z_2^*)$
согласно предложению \ref {pro:2.5}. Так как это решение ограничено,
оно является периодическим в силу условия (a) предложения \ref
{pro:3}, откуда $z^*(\RR^1)=[z_1^*,z_2^*]$ (см.\ шаг 1). Отсюда
$U'(z_i^*)\ne0$ согласно предложению \ref {pro:2.5}. Полученное
противоречие показывает, что $U|_{[c_1,c_2]}=E^*=E_0,$ т.е.\
выполнено условие (b) предложения \ref {pro:3}.

{\it Шаг 4.} Предположим, что выполнено какое-либо из условий (a) и
(b) предложения \ref {pro:3}. Формула (\ref {eq:3.5_0}) доказана на
шаге 1. Докажем (\ref {eq:3.5_1}). Пусть $U''(c_1)=0.$ Тогда для
любого $z\in[c_1-\varepsilon,c_1]$ выполнено $E_0\le U(z)\le
E_0+c(c_1-z)^3$ для некоторого $c>0.$ Отсюда
$\Phi(E)\ge\int_{c_1-((E-E_0)/c)^{1/3}}^{c_1}\frac{dz}{\sqrt{2E-2E_0}}\to+\infty$
при $E\to E_0,$ $E_0<E\le E'.$ Предложение \ref {pro:3} полностью
доказано.
\end{proof}

\begin{propos} [(Локальная обобщенная техническая теорема Бертрана)] \label {pro:4}
Пусть $\Psi=\Psi(z),$ $\rho=\rho(z)$ --- $C^{\infty}$-гладкие функции
на $(a,b).$

{\rm(A)} Пусть $\rho(z_0)\ne0$ для некоторого $z_0\in(a,b).$ Пусть
точка $z_0$ является невырожденным устойчивым положением равновесия
обобщенного уравнения Бертрана при $K=K_0>0,$ т.е.\
$K_0^2=\frac{\Psi(z_0)}{\rho(z_0)}>0$ и
$\rho'(z_0)-\frac{\rho(z_0)}{\Psi(z_0)}\Psi'(z_0)>0.$ Предположим,
что пара $(z_0,K_0)$ удовлетворяет условию {\rm$(\forall)$\loc}\ из
определения \ref {def:rho:closing}. Тогда существуют
$\varepsilon_0,\beta>0,$ такие что функции $\Psi$ и $\rho$ на
$[z_0-\varepsilon_0,z_0+\varepsilon_0]$ удовлетворяют соотношениям
 \begin{equation} \label {eq:3.8_1}
3\Psi''\Psi=4(\Psi')^2, \qquad \Psi\rho'-\Psi'\rho=\beta^2\Psi.
\end{equation}

{\rm(B)} Пусть
выполнены соотношения {\rm(\ref {eq:3.8_1})} на промежутке
$I\subset(a,b).$
Пусть некоторая точка $z_0\in I$ является положением равновесия
обобщенного уравнения Бертрана на $(a,b)$ при некотором $K=K_0\ne0,$
и $\{z_0\}\subsetneq I.$ Тогда в промежутке $I$ выполнены тождества
$\Psi(z)=A_i(z-\zeta)^{1-i^2}$ и $\rho(z) =
\beta^2\frac{(z-\zeta)^4+D}{i^2(z-\zeta)^3},$ где $i\in\{1,2\},$
$A_i,D,\zeta$ --- некоторые константы, такие что $D=0$ при $i=1,$
$z-\zeta\ne0,$ $(z-\zeta)A_1>0,$ $((z-\zeta)^4+D)A_2>0$ при любом
$z\in I.$ При этом любая точка $z\in I$ является невырожденным
устойчивым положением равновесия обобщенного уравнения Бертрана при
$K=\pm\sqrt{\frac{\Psi(z)}{\rho(z)}}$; все ограниченные непостоянные
решения $z=z(\varphi)$ обобщенного уравнения Бертрана на $I$ при
всевозможных значения параметра $K$
являются периодическими с минимальным положительным периодом
$\Phi=2\pi/\beta.$
\end{propos}

\begin{proof} \textit{Шаг 1.}
Обозначим $\C:=\frac2{K^2},$ $\C_0:=\frac2{K_0^2},$
$E_0:=U_{\C_0}(z_0).$ Точка $z_0\in(a,b)$ является невырожденным
устойчивым положением равновесия уравнения
$z''_{\varphi\varphi}(\varphi)=-U'_\C(z)$ при $\C=\C_0$ тогда и только
тогда, когда $U_{\C_0}'(z_0)=0$ и $U_{\C_0}''(z_0)>0.$ Поэтому
$-2U_{\C_0}'(z_0)=R'(z_0)+\C_0W'(z_0)=-2\rho(z_0)+\frac2{K_0^2}\Psi(z_0)=0$
(откуда $\Psi(z_0)\ne0$ ввиду $\rho(z_0)\ne0$) и
$-2U_{\C_0}''(z_0)=R''(z_0)+\C_0W''(z_0)=-2\rho'(z_0)+2\frac{\rho(z_0)}{\Psi(z_0)}\Psi'(z_0)<0.$

Выберем столь малые $\varepsilon,\delta>0,$ удовлетворяющие условию
$(\forall)$\loc\ из определения \ref {def:rho:closing} для пары
$(z_0,K_0),$ что $U_\C''|_{[z_0-\varepsilon,z_0+\varepsilon]}>0$ и
$E_{0,\varepsilon,\C}<E'_{\varepsilon,\C}$ при
$K=\sqrt{2/\C}\in[K_0-\delta,K_0+\delta],$ где обозначено
$E'_{\varepsilon,\C}:=\min\{U_\C(z_0-\varepsilon),U_\C(z_0+\varepsilon)\},$
$E_{0,\varepsilon,\C}:=\min U_\C|_{[z_0-\varepsilon,z_0+\varepsilon]}.$
Для любой пары $(E,\C),$ такой что
 \begin{equation} \label {eq:KE}
K=\sqrt{2/\C}\in[K_0-\delta,K_0+\delta], \qquad
E\in(E_{0,\varepsilon,\C},E'_{\varepsilon,\C}],
 \end{equation}
рассмотрим ограниченное решение $z_{E,\C}(\varphi)$ уравнения
$z''_{\varphi\varphi}=-U_\C'(z)$ с уровнем энергии $E,$ такое что
$z_{E,\C}(\RR^1)\subset[z_0-\varepsilon,z_0+\varepsilon].$ Оно
существует в силу предложения \ref {pro:2.5} и является периодическим
в силу условия $(\forall)$\loc\ из определения \ref
{def:rho:closing}. Запишем по предложению \ref {pro:3} значение
$\Phi(E,\C)$ минимального положительного периода этого решения:
\begin{equation} \label{integreq1}
\Phi(E,\C)=2\int_{z_{1}(E,\C)}^{z_{2}(E,\C)}\frac{dz}{\sqrt{R(z)+\C W(z)+2E}},
\end{equation}
где интегрирование проводится по отрезку между двумя нулями $z_1$ и
$z_2$ знаменателя (а потому, согласно предложению \ref {pro:2.5},
$z_1$ и $z_2$ являются минимумом и максимумом периодического решения
$z_{E,\C}=z_{E,\C}(\varphi)$ уравнения
$z''_{\varphi\varphi}(\varphi)=-U'_\C(z)$).

Используя то, что $z_1=z_1(E,\C)$ и $z_2=z_2(E,\C)$ зануляют
знаменатель в правой части уравнения (\ref{integreq1}), выразим через
них константы $\C$ и $E$ и рассмотрим полученные выражения на
множестве всех пар $(z_1,z_2)\in(a,b),$ таких что $z_0-\varepsilon\le
z_1<z_2\le z_0+\varepsilon$:
\begin{equation} \label{A}
\C(z_1,z_2):=\frac{R(z_2)-R(z_1)}{W(z_1)-W(z_2)}, \quad
2E(z_1,z_2):=\frac{R(z_1)W(z_2)-R(z_2)W(z_1)}{W(z_1)-W(z_2)}.
\end{equation}
Доопределим эти выражения на множестве пар $(z_1,z_1)$ совпадающих
чисел соотношениями
$\C(z_1,z_1):=-\frac{R'(z_1)}{W'(z_1)}=2\frac{\rho(z_1)}{\Psi(z_1)},$
$2E(z_1,z_1):=-R(z_1)+\frac{R'(z_1)}{W'(z_1)}W(z_1).$ Поскольку
$W'(z_0)=\Psi(z_0)\ne0,$ то при достаточно малом $\varepsilon_0>0$
функции (\ref{A})
корректно определены и непрерывны по $(z_1,z_2)$ при
$z_0-\varepsilon_0\le z_1\le z_2\le z_0+\varepsilon_0.$
Следовательно, при достаточно малом $\varepsilon_0>0$ соответствующие
значения параметра $K=\sqrt{\frac2{\C(z_1,z_2)}}$ обобщенного
уравнения Бертрана (для всевозможных пар точек
$z_1,z_2\in[z_0-\varepsilon_0,z_0+\varepsilon_0],$ таких что
$z_1<z_2$) будут принадлежать $\delta$-окрестности числа
$\sqrt{\frac2{\C(z_0,z_0)}}=\sqrt{\frac2{\C_0}}=K_0.$ Отсюда следует,
что для любых $z_1,z_2\in[z_0-\varepsilon_0,z_0+\varepsilon_0],$
таких что $z_1<z_2,$ пара $(E,\C):=(E(z_1,z_2),\C(z_1,z_2))$
удовлетворяет (\ref {eq:KE}) (т.е.\ принадлежит области определения
функции $\Phi=\Phi(E,\C)$ из (\ref {integreq1})), и для любой такой
пары формула (\ref {integreq1}) имеет вид
\begin{equation} \label {eq:PhiEC}
\Phi(E(z_1,z_2),\C(z_1,z_2))=2\int_{z_{1}}^{z_{2}}\frac{dz}{\sqrt{R(z)+\C(z_1,z_2)W(z)+2E(z_1,z_2)}}.
\end{equation}

{\it Шаг 2.} Положим теперь $z_1=c-h,$ $z_2=c+h,$ $z=c+ht,$ где
$c\in(z_0-\varepsilon_0,z_0+\varepsilon_0),$ $0<h\ll1.$ Раскладывая
далее $W$ и $R$ в этих точках в ряды Тейлора по степеням переменной
$h$ и подставляя их в (\ref{eq:PhiEC}) с учетом (\ref{A}),
получаем в пределе при $h\rightarrow0$:
$$
 \lim_{h\to0}\frac{\Phi(E(c-h,c+h),\C(c-h,c+h))}{2\pi}
 =\sqrt{\frac{2W'(c)}{R'(c)W''(c)-W'(c)R''(c)}}
 =:\frac1{\beta(c)}.
$$
Отсюда получаем следующее соотношение в интервале
$(z_0-\varepsilon,z_0+\varepsilon)$:
\begin{equation} \label {eq:second}
R''W'-R'W''=-2\beta^2 W',
\end{equation}
где обозначено $\beta:=\beta(c),$ $R':=R'(c),$ $W':=W'(c),$ $R'':=R''(c),$ $W'':=W''(c).$ Подставляя в
это соотношение $W'=\Psi$ и $R'=-2\rho,$ получаем соотношение
$\Psi\rho'-\Psi'\rho=\beta^2\Psi$ в интервале $(z_0-\varepsilon,z_0+\varepsilon),$ т.е.\ второе
соотношение из (\ref {eq:3.8_1}).

Так как функция (\ref {eq:PhiEC}) является непрерывной в своей
области определения (т.е.\ при $z_0-\varepsilon_0\le z_1<z_2\le
z_0+\varepsilon_0$) и все ее значения попарно соизмеримы ввиду
условия $(\forall)$\loc\ из определения \ref {def:rho:closing}, то
она постоянна и равна своему предельному значению, $2\pi/\beta(z_0).$
В частности, функция $\beta=\beta(c)$ постоянна на интервале
$(z_0-\varepsilon,z_0+\varepsilon).$

\textit{Шаг 3.} Теперь разложим интеграл (\ref{eq:PhiEC}) по степеням
$h$ и рассмотрим коэффициент при $h^{2}$:
\begin{equation} \label {eq:Phi}
\frac{\pi}{4!\beta^3W'}
\left(\frac34(R^{iv}W'-R'W^{iv})+\frac{R'W'''-R'''W'}{W'}
\left(\frac{R'W'''-R'''W'}{8\beta^2}+W''\right)
\right).
\end{equation}
Из равенств (\ref {eq:second}) и $\beta=\const$ получаем:
\begin{equation} \label {eq:iii}
R'W'''-R'''W'=2\beta^2 W''; \qquad
R^{iv}W'-R'W^{iv}=\frac{2\beta^2}{W'} ((W'')^2-2W'W''').
\end{equation}
Используя выражения (\ref {eq:iii}), 
получаем, что при условии $\beta=\const$ коэффициент (\ref {eq:Phi}) при $h^{2}$ в
разложении интеграла (\ref{eq:PhiEC}) зависит только от значений
функции $W$ и ее производных в точке $c$ и равен
$\frac{\pi}{4!}\frac{1}{\beta(W')^2}[4(W'')^2-3W'W'''].$ Приравнивая его к нулю и полагая $W'=\Psi,$
получаем следующее соотношение в интервале
$(z_0-\varepsilon,z_0+\varepsilon)$:
\begin{equation} \label {eq:Psi}
3\Psi''\Psi=4(\Psi')^{2},
\end{equation}
т.е.\ первое соотношение в (\ref {eq:3.8_1}). Часть (A) предложения
\ref {pro:4} доказана.

{\it Шаг 4.} Докажем часть (B) предложения \ref {pro:4} о явном
решении системы дифференциальных уравнений (\ref {eq:3.8_1}). Так как
точка $z_0\in I$ является положением равновесия, то $\Psi(z_0)\ne0.$
Найдем все гладкие решения уравнения (\ref {eq:Psi}) в промежутке
$I,$ такие что $\Psi(z_0)\ne0$: это решения
$\Psi_i(z)=A_i(z-\zeta)^{1-i^2},$ $i=1,2,$
где $A_i \neq 0$ и $\zeta\not\in I$ --- константы. Других гладких
решений, определенных на промежутке $I$ и удовлетворяющих условию
$\Psi(z_0) \neq 0,$ нет, поскольку начальные условия
$(\Psi(z_0),\Psi'(z_0))$ указанных выше решений образуют все
множество $\mathbb{R}^2 \setminus (\{0\} \times \mathbb{R}),$ причем
указанные решения не имеют нулей (даже на всем $\RR$), а потому
являются решениями дифференциального уравнения
$\Psi''=\frac43(\Psi')^2/\Psi$ второго порядка, разрешенного
относительно старшей производной. Наконец, из равенства $\beta =
\const>0$ и второго уравнения в (\ref{eq:3.8_1}), зная
$\Psi=\Psi_i(z),$ находим $\rho$ на промежутке $I$:
$ \rho_{i}(z)={\beta^2}(z-\zeta+D(z-\zeta)^{-3})/{i^2},$ $i=1,2,$
соответственно, где $D$ --- некоторая константа,
такая что $D=0$ при $i=1,$ $D \neq -(z-\zeta)^4$ для любого $z\in I.$
Из равенств $U'_{\C_0}(z_0)=0$ и $\C_0=2/K_0^2>0$ имеем
$\frac{2}{\C_0}=\frac{\Psi(z_0)}{\rho(z_0)}>0,$ откуда знак константы
$A_i$ определяется однозначно по знаку функции $\rho$ на $I,$
$i=1,2.$ Проверка того, что для любой указанной пары функций
$(\Psi_i,\rho_i),$ $i=1,2,$ функция $\Psi_i$ является
$\rho_i$-замыкающей и сильно $\rho_i$-замыкающей на $I,$ причем $\Phi=2\pi/\beta,$ проводится
непосредственно. Предложение \ref {pro:4} доказано.
\end{proof}

\begin{proof}[теоремы \ref{genbert}]
\textit{Шаг 1.}
Предположим сначала, что функция $\Psi$ является полулокально
$\rho$-замыкающей на $(a,b)$ (см.\ определение \ref
{def:rho:closing}). Пусть $z=z(\varphi)$ -- ограниченное решение
обобщенного уравнения Бертрана со значением параметра $K=K_0$ и
уровнем энергии $E'.$ Оно существует ввиду условия $(\exists)$ из
определения \ref {def:rho:closing}. Обозначим $\C_0:=2/K_0^2,$
$a':=\inf z(\RR^1),$ $b':=\sup z(\RR^1),$ $E_0:=\min
U_{\C_0}|_{[a',b']}.$ Тогда $[a',b']\subset(a,b).$ Если
$z_{E,\C_0}(\varphi)$ -- ограниченное решение уравнения при $K=K_0$ с
уровнем энергии $E\in(E_0,E'],$ такое что $z_{E,\C_0}(0)\in(a',b'),$
то в силу предложения \ref {pro:2.5} имеем
$z_{E,\C_0}(\RR^1)\subset[a',b'],$ а в силу условия $(\forall)$\sloc\
из определения \ref {def:rho:closing} решение $z_{E,\C_0}(\varphi)$
периодично. Значит, выполнено условие (a) из предложения \ref {pro:3}
при $U=U_{\C_0},$ а потому, согласно предложению \ref {pro:3},
выполнено условие (b) этого утверждения при $U=U_{\C_0}.$ Мы
утверждаем, что соответствующий отрезок $[c_1,c_2]\subset(a',b')$
является невырожденной точкой локального минимума функции $U_{\C_0}.$
Действительно, в силу (\ref {eq:3.5_1}) и условия $(\forall)$\sloc\
из определения \ref {def:rho:closing}, выполнено $U_{\C_0}''(c_1)>0$
(так как в случае $U_{\C_0}''(c_1)=0$ множество значений функции
периода содержит интервал ввиду (\ref {eq:3.5_1}), а потому не
состоит из попарно соизмеримых чисел, что противоречит условию
$(\forall)$\sloc\ из определения \ref {def:rho:closing}). Значит,
отрезок $[c_1,c_2]\subset(a',b')$ является невырожденной точкой
локального минимума функции $U_{\C_0},$ т.е.\ выполнено условие
$(\exists)$\loc\ из определения \ref {def:rho:closing}. Поэтому
функция $\Psi$ является локально $\rho$-замыкающей на $(a,b).$

Если $\Psi$ является $\rho$-замыкающей или сильно или слабо
$\rho$-замыкающей, то она автоматически является локально
$\rho$-замыкающей (см.\ замечание \ref {rem*}(a)).

{\it Шаг 2.} Пусть теперь функция $\Psi$ является локально
$\rho$-замыкающей на $(a,b),$ и функция $\rho$ не имеет нулей на
$(a,b).$
Пусть $z_0\in(a,b)$ -- невырожденное устойчивое положение равновесия
уравнения при некотором $K=K_0>0$ (оно существует в силу условия
$(\exists)$\loc\ из определения \ref {def:rho:closing}).
Так как выполнено условие $(\forall)$\loc, то, согласно части (A)
локальной обобщенной технической теоремы Бертрана (см.\ предложение
\ref {pro:4}) при $\C=2/K^2$ и $\C_0=2/K_0^2,$ функции $\Psi,\rho$
удовлетворяют системе дифференциальных уравнений (\ref {eq:3.8_1}) в
некоторой окрестности $(z_0-\varepsilon,z_0+\varepsilon)\subset(a,b)$
точки $z_0.$ Пусть $(a_0,b_0)\subseteq(a,b)$ -- максимальный по
включению интервал, содержащий точку $z_0,$ на котором выполнены
дифференциальные соотношения (\ref {eq:3.8_1}).
В силу гладкости функций $\Psi,\rho,$ соотношения (\ref {eq:3.8_1})
выполнены на промежутке $I:=[a_0,b_0]\cap(a,b).$ Отсюда и из
локальной обобщенной технической теоремы Бертрана (см.\ предложение
\ref {pro:4}) следует, что в некоторой окрестности промежутка $I$ в
$(a,b)$ верна одна из формул для пары $(\Psi,\rho),$ указанных в
предложении \ref {pro:4}(B).
Отсюда $I=(a_0,b_0)=(a,b),$ так как в противном случае интервал
$(a_0,b_0)$ немаксимален.

Таким образом, на всем интервале $(a,b)$ верна одна из формул для
пары $(\Psi,\rho),$ указанных в предложении \ref {pro:4}(B). Согласно
предложению \ref {pro:4}(B), $\Psi$ является $\rho$-замыкающей и
сильно $\rho$-замыкающей. Отсюда
$\Psi$ является также полулокально и слабо $\rho$-замыкающей (см.\ замечание \ref {rem*}(a)),
т.е.\ удовлетворяет любому из пяти определений $\rho$-замыкающей функции (см.\ определение \ref {def:rho:closing}).

{\it Шаг 3.} Если $\rho(z)$ является линейной возрастающей на $(a,b)$
функцией, то имеют место обе возможные $\rho$-замыкающие силовые
функции $\Psi_i(z)=A_i(z-\zeta)^{1-i^2},$ $i=1,2,$
причем соответствующие угловые периоды равны
$\Phi_i=2\pi/\beta_i=2\pi/(i\sqrt{\rho'}).$
Если $\rho(z)$ имеет вид $\rho_2(z)$ при $D \neq 0,$ то возможная
$\rho$-замыкающая силовая функция ровно одна: это
$\Psi_2(z)=A(z-\zeta)^{-3},$ причем $\Phi=2\pi/\beta.$ 
В обоих случаях из закона сохранения энергии выводится, что решения $z=z(\varphi)$ соответствующих обобщенных уравнений Бертрана имеют вид, указанный в \S\ref {subsec:orbits}, где $z(r)-\zeta=-\Theta(r)=-\mu^2\theta(r),$ $\Psi(z)=f^2(r)V'(r)=-dV(r(z))/dz,$ $D=\mu^8d.$
Наконец, если
$\rho=\rho(z)$ не принадлежит ни одному из двух указанных видов, то
не будет существовать ни одной $\rho$-замыкающей силовой функции
$\Psi.$

Обобщенная техническая теорема \ref {genbert} Бертрана полностью
доказана.
\end{proof}

\section{Общий случай движения в центральном поле сил} \label {sec:our}

\begin{proof}[теорем \ref {our1} и \ref{our_1_1}]
\textit{Шаг 1.} Опишем рассматриваемую систему. Пусть на поверхности
$S\approx(a,b)\times S^1$ с координатами $(r,\varphi\mod2\pi)$
заданы риманова метрика (\ref {metric}) и потенциал $V=V(r),$ который зависит только от координаты $r.$
Найдём условие на метрику, являющееся критерием возможности
обобщения результата Бертрана на соответствующую систему. Будем
обозначать через $()\dot{}$ производную по $t,$ через $()'_r$ ---
производную по $r,$ а через $()'_{\varphi}$ --- производную по
$\varphi.$

\textit{Шаг 2.} Лагранжиан движения имеет вид $L=\frac{1}{2}\left(\dot{r}^{2}+f^{2}(r)\dot{\varphi}^{2}\right)-V(r),$
а уравнения Эйлера-Лагранжа выглядят так:
\begin{equation}\label{eu1}
ff'_r\dot{\varphi}^{2}-V'_r-\ddot{r}=0, \qquad
(\dot{\varphi}f^{2})\dot{}=0, \qquad K:=\dot{\varphi}f^2=\const.
\end{equation}
При $K=0$ движение происходит по прямой $\{\varphi=\const\}.$ Пусть
далее $K \neq 0.$ Поскольку $\dot{\varphi}=\frac{K}{f^2} \neq 0,$ на
траектории движения $(r(t),\varphi(t))$ можно ввести параметр
$\varphi$ вместо $t.$

\begin{lemma}\label{3.1}
При $K \neq 0$ функция $r=r(\varphi),$ задающая орбиту движения точки
по поверхности $S \approx (a,b)\times S^1$ с метрикой
{\rm(\ref{metric})} в центральном поле с потенциалом $V(r),$
удовлетворяет следующему тождеству:
\begin{equation}\label{imp}
 K^{2}\left(-(\Theta \circ r)_{\varphi\varphi}''+\frac{f_{r}'(r)}{f(r)}\right)=f^{2}(r)V'_r(r),
\end{equation}
где $\Theta = \Theta(r)$ --- произвольная функция, такая что
$d\Theta(r)=\frac{dr}{f^{2}(r)}.$ То есть, функция
$z(\varphi):=-\Theta \circ r(\varphi)$ является решением обобщенного
уравнения Бертрана при $\rho(z)=\frac{f_r'(r(z))}{f(r(z))},
\Psi(z)=f^2(r(z))V'_r(r(z)).$
\end{lemma}

\begin{proof} Производные по $t$ и $\varphi$ связаны соотношениями
$$
\frac{dr}{dt}=\frac{K}{f^{2}(r)}\frac{dr}{d\varphi}=K\frac{d(\Theta\circ r)}{d\varphi}, \qquad
\frac{d^{2}r}{dt^{2}}=\frac{K^{2}}{f^{2}(r)}\frac{d^{2}(\Theta\circ r)}{d\varphi^{2}},
$$
где $d\Theta(r)=\frac{dr}{f^{2}(r)}.$ Первое уравнение в (\ref{eu1})
примет вид
$$
 -(\Theta \circ r)''_{\varphi\varphi}\frac{K^{2}}{f^{2}(r)}+f_{r}'(r)\frac{K^{2}}{f^{3}(r)}
 = V'_r(r).
$$
Или, что эквивалентно, (\ref {imp}),
что и требовалось показать. Лемма доказана.
\end{proof}

В силу леммы \ref{3.1} функции $r=r(\varphi),$ задающие движение
точки по рассматриваемой поверхности в центральном поле сил,
совпадают с решениями уравнений, образующих семейство обобщенных
уравнений Бертрана с параметром $K,$ равным значению кинетического
момента на соответствующих траекториях движения, где $z(r) =
-\Theta(r),$ $d\Theta(r)=\frac{dr}{f^{2}(r)},$ $\rho(z) =
\frac{f'_r(r(z))}{f(r(z))},$
$\Psi(z)=f^2(r)V'_r(r)=-\frac{d}{dz}V(r(z)),$
$\frac{\Psi(z)}{\rho(z)}=\frac{f^3(r(z))}{f'_r(r(z))}V'_r(r(z)),$
эффективный потенциал
$K^2U_{{2}/{K^2}}(z)=\frac{K^2}{2f^2(r(z))}+V(r(z)).$ Отсюда следует,
что каждое из условий $(\exists)$ и $(\exists)$\loc\ на потенциал
$V=V(r)$ (см.\ определение \ref{closing}) равносильно одноименному
условию на функцию $\Psi=\Psi(z)$ (см.\ определение
\ref{def:rho:closing}); а
каждое из условий $(\forall),$ $(\forall)$\loc\ и $(\forall)$\sloc\
на потенциал $V$ равносильно ``рациональному аналогу'' одноименного
условия на функцию $\Psi.$ Здесь рациональный аналог условия на
функцию $\Psi$ получается из этого условия заменой
требования попарной соизмеримости указанных периодов
на (более сильное) требование соизмеримости этих периодов с $2\pi.$
Такие $\rho$-замыкающие функции $\Psi$ назовем \textit{рационально
$\rho$-замыкающими}. По определению \ref {traj}(a)
точка $z_0=-\Theta (r_0)$ является устойчивым положением равновесия
обобщенного уравнения Бертрана при $K=K_0$ тогда и только тогда,
когда окружность $\{r_0\}\times S^1$ является сильно устойчивой
круговой орбитой, такой что значение интеграла кинетического момента
на соответствующей траектории равно $K_0.$ Значит, потенциал $V$ является замыкающим
(соответственно локально, полулокально, сильно или слабо замыкающим) для
метрики {\rm(\ref{metric})} тогда и только тогда, когда функция
$\Psi$ является рационально $\rho$-замыкающей (соответственно
локально, полулокально, сильно или слабо рационально $\rho$-замыкающей). 
В силу обобщенной технической теоремы Бертрана \ref{genbert}, для любой
поверхности $S$ с римановой метрикой (\ref{metric}) любой
полулокально замыкающий потенциал $V$ является локально замыкающим, а
если $f$ не имеет критических точек на $(a,b),$ то все пять классов
замыкающих
потенциалов (см.\ определение \ref {closing}) совпадают и обладают следующими свойствами.

Во-первых, существует не более двух (с точностью до аддитивной и
мультипликативной констант) замыкающих центральных потенциалов
$V(r).$

Во-вторых, существование ровно двух (с точностью до аддитивной и
мультипликативной констант) замыкающих центральных потенциалов $V_1$
и $V_2$ равносильно условию
\begin{equation}\label{xi}
\rho'(z):=\frac{d}{dz}\left(\frac{f'_r(r(z))}{f(r(z))}\right)=\const:=\xi^2>0
\end{equation}
(что эквивалентно условию $f''_{rr}(r)f(r)-(f'_r(r))^{2}=-\xi^{2}$),
где $\xi=\sqrt{\rho'(z)}$ --- положительная рациональная константа,
отвечающая постоянной Бертрана $\beta_i=i\xi$ ($i=1,2,$
в зависимости от вида потенциала) из
теоремы \ref{genbert} (см.\ предложение \ref {pro:4}(B)). В случае
(\ref{xi}) выполнено $z(r) = \frac{f'_r(r)}{\xi^2 f(r)}+\zeta.$
Функция $z(r)$ была определена с точностью до аддитивной константы;
для определенности положим $z(r) := \frac{f'_r(r)}{\xi^2f(r)},$
откуда $\rho(z)=\xi^2z,$ $\Theta(r)=-\frac{f'(r)}{\xi^2f(r)}.$
В силу обобщенной технической теоремы \ref{genbert} Бертрана рационально
$\rho-$замыкающие функции имеют вид $\Psi_i(z)=A_iz^{1-i^2},$
$i=1,2,$ где $A_1z>0,$ $A_2>0.$ Отсюда по лемме \ref {3.1} замыкающие
потенциалы являются гравитационным и осцилляторным, т.е.\ имеют вид
$V_i(r)=(-\int\Psi_i(z)dz)|_{z=-\Theta(r)}=(-1)^i|A_i|\,|\Theta(r)|^{2-i^2}/i+B,$
$i=1,2,$ где $B=\const\in\RR.$

В-третьих, существование ровно одного замыкающего центрального
потенциала равносильно условию $\rho =
-\frac{\Theta}{\mu^{2}}-D\frac{\Theta^{-3}}{\mu^{2}},$ где $D \neq 0,$ $\mu\in\QQ\cap\RR_{>0}$
и выполнено $\rho = \frac{f'_r}{f}.$ Интегрирование равенства
$\frac{f'_r}{f}=-\frac{\Theta+D\Theta^{-3}}{\mu^2}$ по $\Theta$ дает
$\frac{1}{f^2}=\frac{\Theta^2-D\Theta^{-2}+C}{\mu^2},$ где $C$ ---
произвольная константа, откуда
$f(r)={\mu}/{\sqrt{\Theta^2(r)+C-D\Theta^{-2}(r)}}.$ Отсюда рационально
$\rho-$замыкающая функция $\Psi_2(z)=A_2/z^3,$ где $A_2(z^4+D)>0.$
Отсюда замыкающий центральный потенциал
$V_2(r)=(-\int\Psi_2(z)dz)|_{z=-\Theta(r)}=A_2|\Theta(r)|^{-2}/2+B,$ где $B=\const\in\RR.$

Формулы для периодов $\Phi_i$ функций $r=r(\varphi),$ задающих
некруговые неособые ограниченные орбиты, и для кинетических моментов
$K_i,$ $i=1,2,$ для круговых орбит являются повторениями формул из
теоремы \ref{genbert}. Значение кинетического момента $K$ на круговой
орбите $\{r\}\times S^1$ находится из соотношений
$K^2=\frac{\Psi(z)}{\rho(z)}=\frac{f^3(r(z))}{f'_r(r(z))}V'_r(r(z))$ (см.\
теорему \ref{genbert}). Притягивающий центр находится там, где
потенциал минимален. Знак силы равен $-\sgn V_i'(r)=-\sgn\Psi_i(-\Theta(r))
=-\sgn\rho(-\Theta(r))=-\sgn f'(r)=\sgn((\Theta^4(r)+D)\Theta(r)),$ а
потому в притягивающем центре достигаются $\inf f(r)$ и 
$\sup(A_2|\Theta(r)|).$

\textit{Шаг 3.} Итак, мы доказали, что выполнение тождества
$f''f-(f')^2=-\xi^2,$ где $\xi \in \QQ \cap \RR_{>0},$ равносильно
тому, что замыкающих потенциалов ровно два (с точностью до аддитивной
и положительной мультипликативной констант). Предъявим явный вид
функций $f,$ для которых это тождество выполняется.

\begin{lemma}\label{3.3}
Не имеющими нулей решениями уравнения $f''f-(f')^{2}=-\xi^{2}$ при $\xi > 0$ являются следующие функции $f=f(r)$ и только они:
\begin{equation}\label{f1}
\frac{\xi}{\alpha}\sin(\alpha r+\beta), \qquad
\pm\frac{\xi}{\alpha}\sh(\alpha r+\beta), \qquad
\pm \xi r+\beta,
\end{equation}
где $\alpha \neq 0,$ $\beta$ --- произвольные вещественные константы, и
$r$ принадлежит интервалу, в котором $f(r) \neq 0.$
\end{lemma}

\begin{proof} Предположим, что $f'(r) \neq 0$ в окрестности некоторой точки.
Положим $f'=p(f).$ Тогда $f''=p'p.$ Пусть $w=p^{2}.$ Тогда исходное
дифференциальное уравнение примет вид $w'f=2w-2\xi^2.$ Его решением
является функция $w=w(f)=C_{1}f^{2}+\xi^2.$ Отсюда получаем: $\pm
dr=\frac{df}{\sqrt{C_{1}f^{2}+\xi^2}}.$ В зависимости от знака
константы $C_1$ получаются решения (\ref{f1})
при $C_1<0,$ $C_1>0$ и $C_1=0$ соответственно. Начальные условия
$(f(0), f'(0))$ этих решений образуют множество
$\{(\frac{\xi}{\alpha}\sin\beta,\xi\cos\beta),
\pm(\frac{\xi}{\alpha}\sh\beta,\xi\ch\beta), (\beta,\pm\xi)\}
\supset(\mathbb{R}\setminus\{0\})\times\RR,$ где
$\alpha\in\mathbb{R}\setminus\{0\},$ $\beta\in\mathbb{R}.$ Поскольку
$f$ является решением дифференциального уравнения II порядка,
разрешенного относительно старшей производной, с гладкой правой
частью (ввиду $f \neq 0$), то других решений без нулей нет. 
\end{proof}

\begin{remark} \label{rem:3.3}
Аналогичным образом доказывается, что не имеющими нулей решениями
уравнения $-f''f+(f')^2=h=\const \le 0$ являются функции
$f(r)=\frac{\sqrt{-h}}{\alpha}\ch(\alpha r+\beta)$ при $h<0,$
$f(r)=\alpha e^{\beta r}$ при $h=0,$ где $\alpha \neq 0,$ $\beta \in
\mathbb{R}.$
\end{remark}

С учетом замен $\Theta(r)=\mu^2\theta(r),$ $C=\mu^4c,$ $D=\mu^8d$
теоремы \ref{our1} и \ref{our_1_1}, а также формулы орбит из \S\ref
{subsec:orbits} доказаны.
\end{proof}

Мы также доказали теорему \ref {statement}(A),(C), равносильность условий (a) и (b) и единственность набора в (b) из теоремы \ref{statement}(B). Равносильность условий (b) и (c) теоремы \ref{statement}(B) следует из соотношений
$Q_{c_1,d_1}(\theta)=\lambda^{-2}Q_{c,d}(\lambda\theta),$
$r_{c_1,d_1}(\theta)=\lambda(r_{c,d}(\lambda\theta)-r_0),$
$\lambda^{-1}f_{c_1,d_1,k}(\lambda(r-r_0))=f_{c,d,k}(r)$ для
$c:=\lambda^2c_1,$ $d:=\lambda^4d_1,$ $\lambda>0.$ Следствие \ref
{cor:geom}(A),(B),(C) выводится из приводимых в нем формул и
таблицы~\ref {tab:inter}. Следствие \ref {cor:geom}(D) является
известным фактом, относящимся к произвольным поверхностям вращения и
их многомерным аналогам. Первая часть следствия \ref{cor:class}(A)
вытекает из единственности набора $(c,d,k)$ и единственности
изометричного вложения поверхности $S$ в поверхность Бертрана
$(S_{c,d,k},ds^2_{\mu,c,d}),$ см.\ условие (b) теоремы
\ref{statement}(B), ввиду единственности чисел $r_0,\eta$ (при
$\mu>0,$ не обязательно $\mu\in\QQ\cap\RR_{>0}$). Первая часть
следствия \ref{cor:class}(B) вытекает из равносильности условий (b) и
(c) теоремы \ref {statement}(B). Вторые части указанных следствий (о
попарной неизометричности, локальной неизометричности, неподобия,
локального неподобия) вытекают из первых частей и произвольности
множителей $\mu,\lambda>0.$ Сопоставления из следствия
\ref{cor:class}(C) переводят геодезические в геодезические в силу
явных формул для геодезических, см.\ \S\ref {subsec:orbits} при
$A=0.$

\end{fulltext}

\end{document}